\documentclass[12pt]{article}
\usepackage{amsmath}
\usepackage{amssymb}
\usepackage{amsthm}
\newtheorem{theorem}{Theorem}[section]

\newtheorem{example}[theorem]{Example}

\newtheorem{proposition}[theorem]{Proposition}
\newtheorem{remark}[theorem]{Remark}

\def\sd{\partial}
\def\xyb{(\bar x,\bar y)}

\def\dis{\displaystyle}

\def\ran{\rangle}
\def\lan{\langle}

\def\Pi{{\cal P}}

\def\R{I\!\!R}
\def\la{\lambda}
\def\al{\alpha}
\def\del{\delta}
\def\ep{\varepsilon}
\def\sig
\def\sd{\partial}
\def\xb{\overline x}

\def\yb{\overline y}
\def\ub{\overline u}
\def\vb{\overline v}

\def\fb{\overline f}

\def\zb{\overline z}

\def\gr{{\rm Graph}~}

\def\rra{\rightrightarrows}

\def\ep{\varepsilon}
\def\la{\lambda}
\def\La{\Lambda}
\def\al{\alpha}

\def\vf{\varphi}
\def\del{\delta}

\def\ga{\gamma}

\def\sig{\sigma}

\def\epb{\overline{\varepsilon}}

\def\cll{{\mathcal L}}

\def\cb{{\mathcal B}}

\def\cj{{\mathcal J}}
\def\cm{{\mathcal M}}

\def\cu{{\mathcal U}}

\def\ch{{\mathcal H}}

\def\ci{{\mathcal I}}

\def\cl{{\rm cl}}
\def\cone{{\rm cone}}

\def\conv{{\rm conv}~}

\def\ran{\rangle}
\def\lan{\langle}

\def\crit{{\rm Crit}}
\def\dis{\displaystyle}

\title{\bf Towards the theory of strong minimum.\\   A view from variational analysis}

\author{A.D. Ioffe\thanks{Mathematics, the Technion}}
\begin{document}

\maketitle

\noindent{\bf Abstract}.  The key element of the approach  
to the theory of necessary conditions in optimal control discussed in the paper
is reduction of the original constrained problem to unconstrained minimization with subsequent application of 
a suitable mechanism of local analysis to characterize minima  of 
(necessarily nonsmooth) functionals that appear after reduction. 
Using unconstrained minimization at the
crucial step of obtaining necessary conditions  facilitates studies   
of new phenomena and allows to get more transparent and technically simple proofs of known results. In the paper we offer a new proof of the maximum principle for a nonsmooth
optimal control problem (in the standard Pontryagin form) with state constraints
and then prove a new second order condition for a strong minimum in the same
problem but with  data differentiable in control and state variables. The role of variational analysis is twofold. Conceptually, the main considerations behind the reduction
are connected with metric regularity and Ekeland's principle. On the other hand,
the subdifferential calculus offers the main technical instrument for 
proofs of first order conditions.

\vskip 1cm

\section{Introduction}

60 years ago the appearance of  the book by Pontriagin, Boltyanskii, Gamkrelidze and Mischenko \cite{PBGM} 
stimulated a series of studies aimed at finding
general approaches to analysis of necessary conditions in constrained optimization. I just mention two basic ideas that played fundamental role in subsequent developments. 
According to the first, proposed by Dubovitzkii and Milyutin \cite{DM65}, under suitable conditions the cones of
variations of the constraint sets and the cost functional must have empty intersection. The
second idea expressed in the clearest form by Gamkrelidze \cite{RVG65}, was that the image of the solution under a mapping naturally associated with the problem must be a boundary
point of the image of the feasible set under the mapping, and moreover,
again under  suitable conditions, the linearization of the mapping must transform
the cone of feasible variations into a cone not coinciding with the entire range space
of the mapping.   
In both theories  convexity of the cones was a basic requirement and 
final results appeared after application of appropriate separation theorems. 

Further developments revealed limitations of the two approaches and a remarkable difference between  them. The method of
variations of Dubovitzkii-Milyutin proved to be very efficient in the study of higher 
order conditions, in particular in optimal control (see e.g. \cite{FO18,MO,OM,PZ07}). But it does not seem to be suitable for adequate treatment of even first order conditions in problems with nonsmooth data\footnote{Nonconvex subdifferential calculus,
and in particular the ``extremal principle" of Kruger-Mordukhovich that can be viewed as a nonconvex extension of the separation theorem, needs closed sets and lsc functions. But the set of trajectories of a control system is typically not closed in the topology of uniform convergence and its closure contains all relaxed trajectories. Therefore the proofs based on nonconvex separation give maximum principle only for relaxed systems - see e.g. \cite{BM}.}. On the contrary, the boundary point approach (with variational techniques replaced by subdifferential calculus and nonsmooth  controllability criteria) was successfully used to extend
the maximum principle to optimal control problems with nonsmooth data (see e.g. \cite{FHC76,FHC,AI84}) but does not seem to work well with higher order conditions\footnote{The only work known for me where controllability approach is used to get second order conditions is \cite{AMG19}. But it is assumed in the paper that the optimal control takes value in the interior of the set of admissible controls for all $t$.}.

In this paper we shall discuss (mainly in connection with optimal control)
a totally different, non-variational   approach to the study of necessary conditions in constrained optimization. It takes  its origin in the  metric regularity
(or rather metric subregularity) property which is one of the most fundamental
concepts of today's variational analysis. The key element of the approach  is reduction of the original constrained problem to unconstrained minimization with subsequent application of 
a suitable mechanism of local analysis to characterize minima  of 
(necessarily nonsmooth) functionals that appear after reduction.
The possibility to work with unconstrained minimization at the
crucial step of obtaining necessary conditions is the principal advantage of the new theory. 
In an almost obvious way, it opens doors to both the study of second order conditions or to work with first order conditions for nonsmooth problems, sometimes with substantial simplification of arguments. We also hope and believe that the approach will be equally applicable to other problems of the theory of optimal control, not considered in this paper.

 The power of the theory  was already demonstrated 
in several earlier  publications devoted to the maximum principle for optimal control of
systems governed by differential inclusions \cite{AI97a,AI19,RV}. 
Here we consider
the optimal control problem with state constraints and dynamics  in the standard Pontryagin form
$$
\begin{array}{rl}
{\rm minimize} & \ell(x(0),x(T)),\\  
{\rm s.t.} & \dot x = f(t,x,u),\quad u\in U(t),\\
& g(t,x(t))\le 0,\; \Phi (x(0),x(T))\in S.
\end{array}
\leqno ({\bf OC})
$$
The  principal results of the paper include reduction theorems (Theorems \ref{redthm} and \ref{red2}), the maximum principle (Theorem \ref{th1})  and a second order condition for a strong minimum (Theorem \ref{th6}).  Recall that a feasible control process $(\xb(\cdot),\ub(\cdot))$ is called a {\it strong local minimum} in ({\bf OC}) if
for some $\ep>0$  the inequality $\ell(\xb(0),\xb(T))\le\ell (x(0),x(T))$ holds for all feasible pairs $(x(\cdot),u(\cdot))$ satisfying $\| x(t)-\xb(t)\|<\ep$ for all $t$, no matter how far $u(t)$ may be from $\ub(t)$. In this case
$\ub(\cdot)$ is usually called an {\it optimal control} 
and $\xb(t)$  an {\it optimal trajectory} in the problem. 
It should be said however that in proofs
we shall be actually  dealing with a weaker type of minimum
close to what is called {\it Pontryagin minimum}  in \cite{MO}.

The assumptions on the components of the problem differ of course for the first and second order conditions and will be stated in the corresponding parts of the paper. Here we just mention that $x$ and $u$ are elements of some Euclidean spaces, say $\R^n$ and $\R^m$. 
Recall also that
a pair $(x(t),u(t))$, with $u(t)\in U(t)$,  satisfying the differential equation is
called {\it control process} with $u(t)$ being a {\it control function} (or just {\it control}) and $x(t)$ the corresponding {\it trajectory}. Control functions will be always assumed  uniformly bounded  (that is belonging to $L^{\infty}$) which of course does not affect generality of presentation too much. A control process is {\it feasible}
if $x(\cdot)$ satisfies the end point and state constraints.

Reduction theorems to be used here are applied to certain subproblems of ({\bf OC}) obtained by replacement 
of $U(t)$ by smaller and better structured sets. The unconstrained problems that appear as a result
of the reduction theorems resemble the classical Bolza problem with functionals to be minimized looking approximately as follows
$$
\vf(x(\cdot))+ \int_0^T\| \dot x - \psi_0(t,x,u)-\sum_{i=1}^k\al_i\psi_i(t,x,u)\|dt,
$$
where $\al_i$ are nonnegative numbers and $u(t)$ are taken from the mentioned better structured subset of $U(t)$.
The off-integral term $\vf(\cdot)$ in general is a Lipschitz function on the space
of continuous functions. But in the absence of state constraints this is a function of
the end points $(x(0),x(T))$.

Efficiency of the unconstrained reduction technique is then demonstrated by a proof  
of the maximum principle for ({\bf OC}) in Section 4 (Theorem \ref{th1}).
Formally, the theorem is equivalent to the the maximum principle  proved in \cite{RV}. So it could be obtained from the maximum principle for problems with differential inclusions
(as it was done in \cite{RV} and earlier  in \cite{AI97a}
for problems without state constraints). 
But a direct proof based on the reduction theorems  is substantially simpler. 

The final main result, Theorem \ref{th6} proved in the last section, gives a new second order necessary optimality
condition for a strong minimum. 
In subsection 4.3 we discuss in sufficient details the connection of our theorem with two earlier second order conditions for a strong minimum: a very recent result of Frankowska and Osmolovskii
\cite{FO18} and an earlier result of Pales and Zeidan \cite{PZ07}. 
In particular we shall  see that (up to some difference in assumptions with \cite{PZ07})
the condition provided by Theorem \ref{th6} is strictly stronger.  
(Note also that here, as in both papers mentioned above, the optimal control is assumed  only measurable, not piecewise continuous as in 
many other publications dealing with second order conditions.   We  refer to \cite{FO18} also
for a brief account of the literature on second order conditions in optimal control.)


But before we turn to optimal control problems, we introduce in the next section the necessary information from variational analysis including the penalization principles,
central for the unconstrained reduction, and the list of the subdifferential calculus rules used at the last step of the proof of Theorem \ref{th1}.

\vskip 2mm

{\bf Notation}
For $x\in\R^n$ we denote by $\| x\|$ the standard Euclidean norm. 
$C([0,T])$ is the space of continuous functions with the standard norm
$\|x(\cdot)\|_C$ equal to maximum of $\| x(t)\|$ on $[0,T]$. We use the same notation for the space of real-valued functions and for $\R^n$-valued functions and do the same for all other functional spaces introduced below. We hope there will be no confusion
caused and each time the specifics of the space will be clear from the context.  

$W^{1,p}([0,T]$, $1\le p\le\infty$, 
is the Banach space of all absolutely continuous mappings defined on $[0,T]$
with the norm $\|x(\cdot)\|_{1,p}=\| x(0)\| + \|\dot x(\cdot)\|_p$, where $\|\cdot\|_p$
stands for the $L^p$-norm  
with $1\le p\le\infty$.
In what follows we shall use a simpler notation for the space and write just $W^{1,p}$. 

By $\lan\cdot,\cdot\ran$ we denote the inner product in $\R^n$ and the canonical bilinear form in Banach spaces. Again, we hope this will not be a cause of any confusion. We shall also use
the notation $y^*\circ F$ for the composition of a mapping into Banach space $Y$ and the action
of the functional $y^*$. Finally, by $B(x,r)$ we denote the ball of radius $r$ around $x$. The
unit ball aroud the origin will be denoted simply by $B$.

\section{Brief excursion into variational analysis}

Below we offer some information needed for the subsequent discussions. Details can be
found in \cite{AI}.

\vskip 2mm

1.  {\bf   Three basic principles}. We start with the three basic principles. The first is
the well known Ekeland principle which is, by far, one of the most fundamental facts of the variational analysis, in particular the key element in proofs of  many existence theorems. Its statement and proofs can be found in many books, see e.g. \cite{AE,AI,RV}.  

\begin{proposition}[Ekeland's principle]\label{ek}
	Let $X$ be a complete metric space and $f$ a lower semicontinuous function on
	$X$ bounded from below. Let further $f(\xb)\le \inf f+ \ep$. Then for any $\la>0$
	there is a $\ub\in X$ such that $d(\xb,\ub)\le\la$, $f(\ub)\le f(\xb)$ and the function
	$g(x) = f(x) +(\ep/\la)d(x,\ub)$ has a unique global minimum at $\ub$.	
\end{proposition}

The other two principles deal with exact (and necessarily nonsmooth) penalization. 
The first is a simple observation made in 1976 by Clarke (see e.g. \cite{FHC} for the proof).

\begin{proposition}\label{cla}
Let $X$ be a metric space and $\vf(x)$ a function on $X$ which is Lipschitz in a neighborhood of a certain $\xb$. Let further $M\subset X$ contain $\xb$. If $\vf$ attains at $\xb$  a local minimum on $M$, then for any  $K>0$ greater than the Lipschitz constant of $\vf$, the function $\vf(x) + Kd(x,M)$ attains an unconditional local minimum  at $\xb$.  
\end{proposition}

Estimating the distance to the constraint set in an optimization problem may be difficult when
the set is defined by   functional relations. Not surprisingly, Clarke who effectively
used this penalization result to deal with nonsmooth nonlinear programming
problems did  not apply it for optimal control and developed, instead, a totally different
techniques. Later Loewen in \cite{PDL} did apply the proposition  to get maximum principle for a free end point optimal control problem with differential inclusions  but, again, had to use Clarke type techniques for a general problem.

The idea of the following closely connected and more flexible result goes back to  some 1979 papers by the author.  Its proof easily follows from Ekeland's principle 
and Proposition \ref{cla} (see e.g. \cite{AI}).
\begin{proposition}[optimality alternative]\label{alt} Let $X$ be a complete metric space,
	let $f$ be a locally Lipschitz function on $X$ and let $M\subset X$.
	Consider the problem of minimiing $f$ on $M$, and let $\xb$ be a local minimum in the problem.
	Let finally $\vf$ be a nonnegative lsc function on $X$ equal to zero at $\xb$.
	Then the following alternative holds: 
	
	$\bullet$ \ either there is a $\la >0$ such that $\la f + \vf$ has an unconditinal local minimum at $\xb$ (non-singular case);
	
	$\bullet$ \  or there is a sequence of $x_m\not\in \cl M$ converging to $\xb$ and a
	 such that for each $m$ the function
	$\vf(x)+ m^{-1}d(x,x_m)$ attains a global minimum at $X$ (singular case).
	
\end{proposition}

In what follows we  refer to $\vf$ as {\it test function}. The possibility to choose
different test functions adds a lot of flexibility. The price to pay is the necessity to consider a sequence of problems in singular cases but the gain is a substantial extension of the class of problems to be dealt with. For optimal control problems for systems governed by differential inclusions this idea was instrumental  in getting the maximum principle without convexity assumptions on sets of possible velocities  (see \cite{AI97a}).  

\vskip 2mm

2. {\bf Metric regularity}. This is one of the central concepts of variational analysis.
Here we just mention a few facts needed for further discussions.
Let $X$ and $Y$ be metric spaces and $F: X\rra Y$ a set-valued mapping. We  use the same symbol $d(\cdot,\cdot)$ to denote the distance in either  space. It will be always clear  
which space we are talking about. 
Take an  $\xyb\in\gr F$. It is said that $F$ is {\it (metrically) regular} near $\xyb$ if there are $K>0$ and $\ep>0$ such that
$$
d(x,F^{-1}(y))\le K d(y,F(x))
$$
if $d(x,\xb)<\ep$ and $d(y,\yb)<\ep$. It is said that $F$ is 
{\it subregular} at $\xyb$ if the inequality holds with $y=\yb$
and $x\in B(\xb,\ep)$. 
The following is the main (for this paper) example of a regular mapping.

 \begin{proposition}\label{di} Let $S$ be the set of solutions of the differential equation $\dot x = F(t,x)$ defined on $[0,T]$. Assume that we are given an 
$x(\cdot)\in W^{1,1}(\cdot)$, an	$\ep>0$ and a summable $k(t)$ such  that
 $$
 \|F(t,x)-F(t,x')\|\le k(t)\| x-x'\|,\quad{\rm if}\; x,x'\in B( x(t),\ep),\; {\rm a.e.}
 $$
 If 
 $$
 \Big( 1+\int_0^Tk(t)dt\Big)\int_0^T\| \dot{x}(t)-F(t,x(t))\|dt <\ep,
 $$
 then the distance from $x(\cdot)$ to $S$ in $W^{1,1}$ does not exceed
$ K\dis \int_0^T\| \dot{x}(t)-F(t, x(t))\|dt$, where $K$ depends only on $\ep$ and
$k(\cdot)$.
\end{proposition}

 The theorem, probably absent in the literature as stated   
 is an easy consequence of \cite{AI}, Theorem 7.33. Similar results with different estimates
 (and proofs) follow from some earlier publications (e.g. \cite{FHC}, Theorem 3.1.6, , \cite{PDL}, Theorem 2C.5).

\vskip 2mm

3. {\bf Subdifferentials}. There are several types of subdifferentials used in local variational analysis. We shall basically work with   the $G$-subdifferential 
(coinciding with the limiting Fr\`echet subdifferential in finite dimensional spaces) and Clarke's generalized
gradient. In what follows the symbol $\sd$ will be used for the first and $\sd_C$ for the second.
These are the only ``good" subdifferentials that make sense and work in arbitrary Banach spaces. Moreover,
if $X$ is a separable Banach space then the $G$-subdifferential is  the minimal subdifferential 
 with the following properties (among others)

\vskip 1mm

$\bullet$ if $f$ is lsc and attains a local minimum at $x$, then $0\in\sd f(x)$; 

\vskip 1mm

$\bullet$ if $f$ is locally Lipschitz then $\sd f(x)\neq\emptyset$ and the mapping
$x\to\sd f(x)$ is bounded-valued and norm-to-weak$^*$ usc; 

\vskip 1mm

$\bullet$ if $f$ is Lipschitz near $x$ then  
$\conv\sd f(x)$ coincides with Clarke's generalized gradient  $\sd_Cf(x)$ of $f$ at $x$;

\vskip 1mm

$\bullet$ if $f$ is convex, then $\sd f(x)$ coincides with the subdifferential in the sense of
convex analysis; if $f$ is strictly differentiable at $x$, then $\sd f(x) = \{f'(x)\}$.

\vskip 1mm


If $S\subset X$ is a closed set and $x\in S$, then
$N(S,x)=\cone ~\sd d(\cdot,S)(x)$  is  the {\it normal cone} to $S$ at $x$.

Here are some  calculus rules for  $G$-subdifferentials to be used in the paper:

\vskip 1mm

\vskip 1mm 

$\bullet$ \ $\sd (\la \sd f(x))= \la f(x)$ ($\la >0 $);

\vskip 1mm

$\bullet$ \ If $f= f_1+f_2$, where both function are lsc and at least one of them is  Lipschitz near $x$, then
$\sd f(x) \subset \sd f_1(x)+\sd f_2(x)$.

\vskip 1mm

$\bullet$ \ If $f(x,y)= f_1(x)\cdot f_2(y)$ and both $f_1$ and $f_2$ are nonnegative and Lipschitz near $\xb$ and $\yb$ respectively then, $\sd f\xyb= f_1(\xb)(\{0\}\times \sd f_2(\yb)) + f_2(\yb)(\sd f_1(\xb)\times\{0\})$.

\vskip 1mm

$\bullet$ If $X$ is a closed subspace  of $L^{\infty}([0,T],\R^n)$, $f(t,x)$ is measurable in $t$
and $k(t)$-Lipschitz in $x$ in the $\ep$-neighborhood  of $\xb(t)$ a.e. (with summable $k(\cdot)$) and $f(x(\cdot))=\int_0^Tf(t,x(t))dt$, then $\sd f(x(\cdot))\subset\int_0^T\sd_Cf(t,x(t))dt$ in the sense that
 for any $x^*\in\sd f(\xb(\cdot))$ there is a summable  $\xi(t)\in\sd_Cf(t,x(t))$ a.e. such that for all $h(\cdot)\in X$
$$\lan x^*,h(\cdot)\ran =\int_0^T\lan \xi(t),h(t)\ran dt.$$  

\vskip 1mm

$\bullet$ \ If $f=g\circ F$, $F: \R^n\to\R^m$ and both $g$ and $F$  are Lipschitz near $y=F(x)$, and $x$ respectively, then
$$
\sd f(x) \subset \bigcup_{y^*\in\sd g(y)}(y^*\circ F)(x).
$$

\vskip 1mm

$\bullet$ \ If $f(x)= \max_if_i(x)$, $i=1,\ldots,k$,  and all $f_i$ are Lipschitz near $x$, then
$$
\sd f(x) \subset\big\{\sum_{i\in I(x)}\al_i x_i^*:\; x_i^*\in\sd f_i(x),\; \al_i\ge 0,\; \sum\al_i=1    \big\},
$$ 
where $I(x)=\{i: \; f_i(x)=f(x)\}$.

\vskip 1mm

$\bullet$ \ If $X= C[0,T]$ and $f(x(\cdot))=\max_{t\in[0,T]}g(t,x(t))$, where $g$ is an usc real-valued function
satisfying $|g(t,x)-g(t,x')|\le k(t)\| x-x'\|$ if $\| x- x(t)\|<\ep$ with summable $k(\cdot)$ and $\ep>0$,  then $\sd f(x(\cdot))$ consists of measures
$\nu$ such that $\nu(dt)=\ga(t)\mu(dt)$, where $\mu$ is a probability measure
supported on $\Delta=\{t:\; g(t,x(t))=f(x(\cdot)) \}$ and $\ga(\cdot)$ is a measurable
selection of the set-valued mapping $t\to \bar{\sd} g(t,x(\cdot))(x(t))$, where
$$
\bar{\sd}g(t,x)=\{\lim y_m:\; y_m\in\sd g(t_m,\cdot)(x_m),\; t_m\to t,\; x_m\to x,\; g(t_m,x_m)\to g(t,x)    \}.
$$ 
We refer to \cite{FHC,AI,JPP,RW} for further details.

3. {\bf Tangent cones}. Let again $X$ be a Banach space,  $Q\subset X$ is closed, $x\in Q$.
The {\it (classical) tangent cone} $T(Q,x)$ to $Q$ at $x$ is the collection of all
$h\in X$ such that $d(x+th,Q)= o(t)$.

If $h\in T(Q,x)$ then the collection $T^2(Q,x;h)$  of $v\in X$ such that  $d(x+th+t^2v, Q)=o(t^2)$
 is the {\it second order tangent set} to $Q$ at $x$ along $h$.
We denote by $T_0(Q,x)$ the collection of
$h\in T(Q,x)$ for which $T^2(Q,x;h)\neq\emptyset$.

\vskip 2mm

4. {\bf Measurability}. Recall that a set $Q\subset [0,T]\times \R^n$ is {\it $\cll\times\cb$-measurable} if it belongs to the $\sigma$-algebra generated by all
products $\Delta\times V$, where $\Delta\subset [0,T]$ is Lebesgue measurable
and $V\subset \R^n$ is open. A set-valued mapping $V(t)$ from $[0,T]$ into $\R^n$
is {\it $\cll\times\cb$-measurable} (or just {\it measurable}) if its graph is $\cll\times\cb$-measurable. A single-valued measurable mapping $v(t)$ from $[0,T]$ into $\R^n$ is a {\it measurable selection} of $V(\cdot)$ if $v(t)\in V(t)$ for almost every $t$. For us the most important properties of measurable set-valued mappings are:

$\bullet$ if $V(t)$ is measurable, then $\cl V(t)$ (where $cl V$ is the closure of $V$)
is also measurable;

$\bullet$ a measurable set-valued mapping has a measurable selection (Aumann selection theorem).

Specifically, we shall need the following consequence of the last property.

\begin{proposition}\label{mea}
Let $U(t)$ be an $\cll\times\cb$-measurable set-valued mapping from $[0,T]$ into $\R^m$
and $f(t,u)$ an $\cll\times\cb$-measurable extended-real valued function on $[0,t]\times\R^m$. Assume further that a measurable selection $\ub(\cdot)$ of $U(\cdot)$ be
given such that $f(t,u(t))$ is a summable function and
$$
\int_0^Tf(t,\ub(t))dt \ge \int_0^Tf(t,u(t))dt
$$
for any measurable selection $u(\cdot)$ of $U(\cdot)$ for which the integral in the right side of the inequality makes sense. Then for almost every $t$ the inequality $f(t,\ub(t))\ge f(t,u)$ holds for all $u\in U(t)$.

\proof Assume the contrary, and let $V(t)=\{u\in U(t): f(t,u)>f(t,\ub(t))\}$.  Then the graph of $V(\cdot)$ is $\ell\times\cb$-measurable and $V(t)\neq\emptyset$ on a set $\Delta$ of positive measure. Let $v(t)$ be a measurable selection of $V(\cdot)$ defined on $\Delta$,
and let $u(t)$ coincides with $v(t)$ on $\Delta$ and with $\ub(t)$ on $[0,T]\backslash\Delta$. Then $f(t,u(t))\ge f(t,\ub(t))$ for all $t$, hence the integral 
of $f(t,u(t))$ makes sense and we come to a contradiction.\endproof 	
\end{proposition}

\noindent All other results relating to measurable set-valued mappings to be used in the paper
are immediate consequences of the definitions. We refer to \cite{CV,RV} for more details.

\section{Reduction to unconstrained problems}

 Let $(\xb(\cdot),\ub(\cdot))$ be a strong local minimum in ({\bf OC}). This means that
there is a $\ep>0$ such that $\ell(\xb(0),\xb(T))\le\ell(x(0),x(T))$ for any
feasible $(x(\cdot),u(\cdot))$ such that $\| x(t)-\xb(t)\|<\ep$ for all $t$. 
The following
assumptions on components of ({\bf OC}) will be adopted throughout the paper (and either strengthened or supplemented with additional assumptions whenever necessary):   

\vskip 1mm

(H$_1$)  $\ell$ is  a real-valued function, Lipschitz  in a neighborhood of $(\xb(0),\xb(T))$;

\vskip 1mm 

(H$_2$)  there are $\ep_0>0$, $\del_0\ge 0$ and a summable $k_0(t)$ such that 
for almost every $t$ the relations
$$
\|f(t,\xb(t),u)\|\le k_0(t);\quad
\| f(t,x,u)-f(t,x',u)\|\le k_0(t)\| x-x'\|\quad {\rm a.e.}
$$
hold for  $x,x'\in B(\xb(t),\ep_0)$, $u\in U_0(t)= B(\ub(t),\del_0)\cap U(t)$;

\vskip 1mm

(H$_3$)  the set-valued mapping $t\to \{(u,f(t,\xb(t),u)):\; u\in U(t)\}$ is  $\cll\times\cb$-measurable;

\vskip 1mm

(H$_4$) $g(t,x)$ is upper semicontinuous in both variables and there is a $\rho_0>0$ such that $g(t,\cdot)$ is $\rho_0$-Lipschitz  on $B(\xb(t),\ep_0)$ for every $t\in [0,T]$;

\vskip 1mm

(H$_5$) $S\subset \R^r$ is a closed set,   $\Phi: \R^n\times\R^n\to\R^r$ is continuously differentiable near $(\xb(0),\xb(T))$ 
and $\Phi'(\xb(0),\xb(T))$ is a linear operator onto $\R^r$.    .

\vskip 1mm

\noindent (We do not exclude the possibility that $\del_0=0$. In fact, this is exactly the case we shall be dealing with in the proof of the maximum principle in the next section.)

The reduction theorems we are going to state and prove in this section actually apply
not to ({\bf OC}) but to its subproblems defined as follows. Denote by $\cu$ the collection
of all mesurable selections $u(\cdot)$ of $U(\cdot)$ for which there are $\ga>0$ and
summable $k(t)$ such that for almost every $t$ the inequalities
$$
\| f(t,\xb(t),u(t))\|\le k(t),\quad \|f(t,x,u(t))-f(t,x',u(t))\|\le k(t)\| x-x'\|
$$
hold for all $x,x'\in B(\xb(t),\ga)$. By (H$_2$) $\ub(\cdot)\in\cu$. We shall always assume that
$\ga\le\ep_0$ and $k(t)\ge k_0(t)$ a.e..
It will  be also convenient to occasionally denote the $\ga$ and $k(\cdot)$ associated with a given $u(\cdot)\in\cu$ by
$\ga_u$ and $k_u(\cdot)$.

Take now a finite collection
$\{u_1(\cdot),\ldots,u_k(\cdot)\}$ of elements of $\cu$, and let
$$
U_k(t) = U_0(t)\cup\{u_1(t),\ldots,u_k(t)  \}
$$
(with $U_0(\cdot)$ from (H$_2$)).  The subproblem we shall work with is
$$
\begin{array}{rl}
{\rm minimize} & \ell(x(0),x(T)),\\
{\rm s.t.} & \dot x = f(t,x,u),\quad u\in U_k(t),\\
& g(t,x(t))\le 0,\; \Phi (x(0),x(T))\in S
\end{array}
\leqno ({\bf OC}_k)
$$
Clearly $(\xb(\cdot),\ub(\cdot))$ is a strong local minimum in ({\bf OC}$_k$).

In what follows we denote by $\cu_k$  the collection of   measurable $u(\cdot)$ such that
$u(t)\in U_k(t)$ a.e..  It is clear that
$\cu_k\subset \cu$. Let further
$X_k$ denote  the collection of all pairs $(x(\cdot),u(\cdot))$ with $\| x(\cdot)-\xb(\cdot)\|_C\le\ep_0$
and $u(\cdot)\in \cu_k$  satisfying the equation $\dot x=f(t,x,u)$.
Thus ({\bf OC}$_{k}$) is the problem of minimizing  $\ell(x(0),x(T))$ on the set
$M_k$ of elements of $X_k$ satisfying $\Phi(x(0),x(T))\in S$. 
We shall endow $X_k$   
with the $C([0,T])\times L^1$-metric. As $U_{k}(t)$ are closed sets bounded by  summable functions, $X_k$ is a complete metric space.

We shall also consider the space $Z_k$ of $(k+2)$-tuples 
$z=(x(\cdot),u(\cdot),\al_1,\ldots,\al_k)$ with $x(\cdot)\in W^{1,1}$,
$u(\cdot)\in\cu_k$,    $\al_i\ge 0$ and $\sum\al_i\le 1$. 
Unlike $X_k$, the $u(\cdot)$-components of elements of $Z_k$ will be considered with
the $L^{\infty}$-topology of uniform convergence almost everywhere. 
Set
$$
\psi(x(\cdot))= \max\{\ell(x(0),x(T))-\ell (\xb(0),\xb(T)),\max_{0\le t\le T}g(t,x(t)) \}; 
$$
and 
$$
\begin{array}{l}
J_k(z)=J_k(x(\cdot),u(\cdot),\al_1,\ldots,\al_k)\\
\quad\quad\quad= \dis\int_0^T\|\dot x(t)-f(t,x(t),u(t))-\dis\sum_{i=1}^k\al_i(f(t,x(t),u_i(t))-f(t,x(t),u(t)))\|dt.
\end{array}
$$
We are ready to state and prove the first reduction theorem.
\begin{theorem}\label{redthm}
We posit (H$_1$)-(H$_5$). If $(\xb(\cdot),\ub(\cdot))$ is a strong local minimum in 
({\bf OC}$_k$), then the following alternative holds with some sufficiently big $K>0$:

-- either there is a $\la>0$ such that the functional  
$$
\cj_{k0}(z)= \la\psi (x(\cdot)) + d(\Phi(x(0),x(T)),S) +J_{k}(x(\cdot),u(\cdot),\al_1,\ldots,\al_k)
$$
attains a local minimum on $Z_k$ at $\zb=(\xb(\cdot),\ub(\cdot),0,\ldots,0)$;

--or there is a $K>0$ and a sequence $(\xb_m(\cdot),\ub_m(\cdot))\subset X_k\backslash M_k$, $m=1,2,\ldots$ converging to $(\xb(\cdot),\ub(\cdot))$ and such that for any (sufficiently large) $m$ the functional
$$
\begin{array}{l}
\cj_{km}(z)= d(\Phi(x(0),x(T)),S) +KJ_{k}(x(\cdot),u(\cdot),\al_1,\ldots,\al_k)\\
\qquad\qquad+ \  m^{-1}\Big(\|x(\cdot) -\xb(\cdot)\|_C +\dis\int_0^T\big(\| u(t)-\ub_m(t)\|   
+\sum\al_i\|u_i(t)-\ub_m(t)\|       \big)dt    \Big)
\end{array}
$$
attains a local minimum on $Z_k$ at $z_m=(\xb_m(\cdot),\ub_m(\cdot),0,\ldots,0)$.
\end{theorem}

\proof
1.  We start with an almost obvious remark that ({\bf OC}$_k$) can be equivalently
reformulated as the problem of minimuzing $\psi$ on $X_k$ subject to $\Phi(x(\cdot),u(\cdot))\in S$:  
$$
{\rm minimize}\quad \psi(x(\cdot)),\quad{\rm s.t.}\; \dot x=f(t,x,u),\; u\in U_{k}(t),\;
\Phi(x(0),x(T))\in S.
$$
We shall continue to refer to this last formulation of the problem as ({\bf OC}$_{k}$).  
 We can also assume without loss of generality that $\ell(\xb(0),\xb(T)=0$. 
 
Next we apply the optimality alternative to the problem.  
(Note that $\psi(\cdot)$ is  Lipschitz in a neighborhood of
$\xb(\cdot)$.)   Then either there is a $\la_0>0$ such that $(\xb(\cdot),\ub(\cdot))$ is a local minimum of 
$$
\ci_0(x(\cdot))=\la_0\psi(x(\cdot))+ d(\Phi(x(0),x(T)),S)
$$
on $X_k$ (that is subject to $\dot x=f(t,x,u),\; u\in U_{k}(t)$)
({\it non-singular case}) or there is a sequence of pairs $(\xb_m(\cdot),\ub_m(\cdot))\in X_k\backslash M_k$ converging
to $(\xb(\cdot),\ub(\cdot))$, and such that
$$
\ci_m(x(\cdot),u(\cdot))= d(\Phi(x(0),x(T)),S) + m^{-1}\Big(\|x(\cdot)-\xb_m(\cdot)\|_C+ \dis\int_0^T\| u(t)-\ub_m(t)\|dt\Big)
\qquad\qquad\qquad\qquad\qquad\qquad\qquad\qquad\qquad\qquad\qquad
\ge d(\Phi(\xb_m(0),\xb_m(T)),S) >0
$$  
attains a global minimum on $X_k$ at  $(\xb_m(\cdot),\ub_m(\cdot))$ ({\it singular case}). 


\vskip 1mm

2. Fix a $u(\cdot)\in \cu_k$  and $\al_i\ge 0$ with $\sum\al_i<1$, and set
$\ga=\min\{\ga_{u_1},\ldots,\ga_{u_k}\}$ and $k(t)=\max\{k_{u_1}(t),\ldots,k_{u_k}(t)  \}$.
Let further $x(\cdot)$ be a solution of the  following differential equation:
$$
\dot x= f(t,x,u(t))+\sum_{i=1}^k\al_i(f(t,x,u_{i}(t))-f(t,x,u(t))).
\eqno (1)
$$
 If $x(\cdot)$ is sufficiently close to $\xb(\cdot)$, say  $\| x(\cdot)-\xb(\cdot\|_C<\ga/2$,
 then there is a sequence $(x_s(\cdot),u_s(\cdot))$ 
satisfying the original equation $\dot x=f(t,x,u)$ with $x_s(\cdot)$ uniformly converging to $x(\cdot)$ and $u_s(t)\in\{u(t),u_1(t),\ldots,u_k(t)\}$ almost everywhere. 

This is immediate from a number of relaxation theorems (see e.g \cite{PDL}, Theorem 2F.2
or \cite{RV}, Theorem 2.7.2). The absence of end point constraints makes construction of
the desired sequence especially easy. 
All we need is to break $[0,T]$ into $s$ equal intervals $\Delta_j$, then to choose  in each of them  $k$ disjoint
subsets $\Delta_{ij}$ with measures respectively equal to $\al_iT/s$ and set
$$
u_s(t) = \big(1-\sum_{i=1}^k\al_{is}(t)  \big)u(t) + \sum_{i=1}^k\al_{is}u_i(t),
$$
where
$$
\al_{is}(t) = \left\{\begin{array}{l}1,\ {\rm if}\; t\in\cup_{j}\Delta_{ij};\\ 0,\ {\rm otherwise}. \end{array}\right.
$$
Then for any $i=1,\ldots,k$ the sequence $(\al_{is}(\cdot))$ weakly (e.g. in $L^1$) converges to the function identically equal to $\al_i$ and it is not a difficult matter to deduce 
(taking (H$_2$) into account) that the sequence $(x_s(\cdot))$ of solutions of equations
$\dot x=f(t,x,u_s(t))$ with $x_s(0)=x(0)$ uniformly converges to $x(\cdot)$.

As $\ci_m(x_s(\cdot),u_s(\cdot))\ge \ci_m(\xb_m(\cdot,\ub_m(\cdot)))$,
it follows that in the non-singular case $(\xb(\cdot),\ub(\cdot),0,\ldots,0)$ is a local minimum in the problem of minimizing  $\ci_0$ on the set of $z=(x(\cdot),u(\cdot),\al_1,\ldots,\al_k)\in Z_k$ satisfying (1).  
For the same reason in the singular case the obvious inequality
$$
\| u_s(t)-\ub_m(t)\|\le (1-\sum_{i=1}^m\al_{is}(t))\| u(t)-\ub_m(t)\|+\sum_{i=1}^m\al_{is}(t)\| u_i(t)-\ub_m(t)\|
$$
leads to the conclusion that $(\xb_m(\cdot),\ub(\cdot),0,\ldots,0)$ is a global minimum in the problem of minimizing the functional
$$
\begin{array}{l}
d(\Phi(x(0),x(T)),S) + m^{-1}\Big(\|x(\cdot) -\xb_m(\cdot)\|_C\\
\qquad\qquad\qquad
 +\dis\int_0^T\big( \| u(t)-\ub_m(t)\|  
 +\sum\al_i\|u_i(t)-\ub_m(t)\|       \big)dt    \Big)
\end{array}
$$ 
on the same set of $z=(x(\cdot),u(\cdot),\al_1,\ldots,\al_k)\in Z_k$ satisfying (1).  
If $\al_1=\ldots=\al_k=0$ this quantity coincides with $\ci_m(x(\cdot),u(\cdot))$, so it can be viewed as an extension of $\ci_m$ to $Z_k$ and we can continue to keep the notation
$\cm_m$ for this functional.

\vskip 1mm

3. Given a $(u(\cdot),\al_1,\ldots,\al_k)\in \cu_k\times \R_+^k$ with $\sum\al_i<1$,  we denote by\\  $Q(u(\cdot),\al_1,\ldots,\al_k)\subset W^{1,1}$ the set of solutions of (1).  
As follows from Proposition \ref{di} (in view of (H$_2$)), there are $K_0>0$, $\ga>0$ such that the inequality
$$
d(x(\cdot),Q(u(\cdot),\al_1,\ldots,\al_k))\le K_0 J_k(x(\cdot),u(\cdot),\al_1,\ldots,\al_k)
$$
holds for $(x(\cdot),u(\cdot),\al_1,\ldots,\al_k)\in Z_k$ satisfying
$$
\|x(\cdot)-\xb(\cdot)\|_C+ \big(1+\int_0^Tk(t)dt\big)J_k(x(\cdot),u(\cdot),\al_1,\ldots,\al_k)<\ga.
$$

It remains to apply
 Proposition \ref{cla} to the result obtained at the end of the previous  step of the proof.
In the non-singular case it follows that $(\xb(\cdot),\ub(\cdot),0,\ldots,0)$ is a local minimum of $\ci_0(x(\cdot),u(\cdot))+ KJ_k(x(\cdot),u(\cdot),\al_1,\ldots,\al_k)$ on $Z_k$
for some $K>0$. 
This implies the ``either" part of the statement with $\la=\la_0$ if $K\le 1$ and $\la=\la_0/K$ if $K>1$.

 Likewise, in the singular case we deduce that $(\xb_m(\cdot),\ub_m(\cdot),0,\ldots,0)$ is a local minimum of $\ci_m(x(\cdot),u(\cdot),\al_1,\ldots,\al_k)$ and the ``or" part of the statement follows as well.  \endproof

The conclusion of the theorem in the non-singular case may be not fully satisfactory 
in certain cases for the following reason. The necessary optimality condition
for $\cj_k$ obtained by application of the subdifferential calculus to $\cj_{k0}$ will inevitably include the limiting or Clarke's subdifferential of $\psi$ at $\xb(\cdot)$.  But this subdifferentials contain vectors obtained from the analysis of behavior of the function
at points (close to $\xb(\cdot)$) at which the value of $\psi$ is strictly smaller than
$\psi(\xb(\cdot))$. If such points do exist, they lie are outside of the feasible domain of the problem, so taking
them into account may only decrease precision of the necessary condition. 

To avoid such a possibility, we shall slightly modify the functional obtained  in the non-singular case at the first part of the proof of the theorem and consider the sequence of problems of
minimizing $\la_0\psi_m(x(\cdot)) + d(\Phi(x(0),x(T)),S)$ on $X_k$, where
$$
\psi_m(x(\cdot))= \max\{\ell(x(0),x(T))+m^{-2},\max_{0\le \la \le T}g(t,x(t))  \}
$$
Obviously $\xb(\cdot)$ may no longer be  a local minimum in this problem. But the value of the
new functional at $\xb(\cdot)$ can exceed the minimal value in the problem at most by $m^{-2}$. 

So by  Ekeland's principle,
for any suficiently  big $m$ there is a pair $(\xb_m(\cdot),\ub_m
(\cdot))\in X_k$ such that $\| \xb_m(\cdot)-\xb(\cdot)\|_C<m^{-1}$, $\|u(\cdot) -\ub_m(\cdot)\|_{L^1}\le m^{-1}$ and the function
$$
\la_0\psi_m(x(\cdot)) + d(\Phi(x(0),x(T)),S)+ m^{-1}(\| x(\cdot)-\xb_m(\cdot)\|_C
+\int_0^T\|u(t)-\ub_m(t)\|dt)
$$
attains a local minimum on $X_k$ at $(\xb_m(\cdot),\ub_m(\cdot) )$.

Thus we arrive to a series of problems that are very similar to what we have in the singular  cases. The subsequent arguments in the proof of the theorem can be applied to these problems as well, and we arrive at the following conclusion.

\begin{theorem}\label{red2}
	We posit (H$_1$)-(H$_5$). If $(\xb(\cdot),\ub(\cdot))$ is a strong local minimum in 
	({\bf OC}), then  there are $\la\ge 0$,  $K>0$ and a sequence $(\xb_m(\cdot),\ub_m(\cdot))\subset X_k$, $m=1,2,\ldots$ converging to $(\xb(\cdot),\ub(\cdot))$ and such that the functional
	$$
	\begin{array}{l}
	\la\psi_m(x(\cdot)+ d(\Phi(x(0),x(T)),S)
	+KJ_{k}(x(\cdot),u(\cdot),\al_1,\ldots,\al_k)\\
	\qquad\quad
	+ \  m^{-1}\Big(\|x(\cdot) -\xb_m(\cdot)\|_C +\dis\int_0^T\big(\| u(t)-\ub_m(t)\|dt + \sum\al_i\| u_i(t)-\ub_m(t)\|\big)dt         \Big)
	\end{array}
	$$
	attains a local minimum on $Z_k$ at $z_m=(\xb_m(\cdot),\ub_m(\cdot),0,\ldots,0)$
	and  either $\la>0$ or $\la = 0$ and 
	$\Phi(\xb_m(\cdot),\ub_m(\cdot))\not\in S$.
\end{theorem}

\section{First order condition (maximum principle)} 

\subsection{Basic unconstrained model}
By that we mean the problem of minimization the functional
$$
J(x(\cdot),\al_1,\ldots,\al_k)= \vf(x(\cdot))
+\int_0^T\|\dot x(t)-\psi_0(t,x(t))-\sum_{i=1}^k\al_i\psi_i(t,x(t))\|dt
$$
on $W^{1,1}\times \R_+^k$,
where $\vf$ is a Lipschitz function on $C[0,T]$ and $\psi_i: [0,T]\times \R^n\to \R^n$
are Carath\`eodory functions satisfying 
$$
\| \psi_i(t,x)\|\le R(t),\quad \|\psi_i(t,x)-\psi_i(t,x')\|\le R(t)\| x-x'\|,\quad{\rm if}\; x,x'\in B(\xb(t),\epb)
$$
 with some $\epb>0$ and summable $R(t)$. 
 
If $\vf$ is a function of the end points $(x(0),x(T))$, this is a very specific case of 
the well studied ``generalized Bolza problem" (see e.g. \cite{FHC,IR96,PDL}). 
The presence of nonnegative parameters $\al_i$ does not add much difficulty to analysis. A case of a
general $\vf$ was recently considered  in \cite{AI19}, and the first order necessary
condition for a minimum of $J$ proved in this section is immediate from the main result of \cite{AI19}. But we prefer
to give a separate proof here which, as we have already mentioned, is noticeably simpler.

So let $(\xb(\cdot),(0,\ldots,0))$ be a local minimum of $J$ in $W^{1,1}\times \R_+^k$. 
By the assumption $\| \psi_i(t,x)\|\le R(t)$
for all $x\in \xb(t)+\epb B$ a.e. on $[0,T]$. It follows 
(from Proposition \ref{cla}) that for a sufficiently large $K>0$  the functional
$$
I(x(\cdot),\al_1,\ldots,\al_k)= J(x(\cdot),\al_1,\ldots,\al_k) + K\sum_{i=1}^k\al_i^-
$$
attains a local minimum on $W^{1,1}\times \R^k$
at $(\xb(\cdot),0,\ldots,0)$. 
Here $\al^-=\max\{0,-\al\}$.  Thus $0\in\sd I(\xb(\cdot),0\ldots,0)$.

Applying the calculus rules for the $G$-subdifferential,
 we find a measure $\nu\in\sd\vf(\xb(\cdot))$,  measurable 
$\R^n$-valued functions $p(t)$  satisfying $\|p(t)\|\le 1$ and measurable
$q(t) \in \sd_C\lan p(t), \psi_0(t,\cdot)\ran(\xb(t)),\; {\rm a.e.}$ such that
$$
\int_0^T\big(h(t)\nu(dt)+( \lan \dot h(t),p(t)\ran - \lan h(t),q(t)\ran)dt \big)=0,
\eqno (2)
$$
for all $h(\cdot)\in W^{1,1}$, and \ $0\in \int_0^T\lan p(t),\psi_i(t,\xb(t)\ran dt)+[-K,0]$, that is
$$
\int_0^T\lan p(t),\psi_i(t,\xb(t))\ran  dt\le 0,\quad, i=1,\ldots,k.
\eqno (3)
$$
 Setting $h(t)= h(0)+\int_0^t \dot h(s)ds$, we get
$$
\begin{array}{l}
\dis\int_0^Th(t)\nu(dt)= \lan h(0),\nu(\{0\})\ran + \lan h(T),\nu(\{T\})\ran +\dis\int_0^Th(t)\tilde{\nu}(dt)\\ 
\qquad\qquad=\big{\lan} h(0),\nu(\{0\})+\dis\int_0^T\tilde{\nu}(dt)\big{\ran} + 
 \lan h(T),\nu(\{T\})\ran +\dis\int_0^T\big{\lan}\dot h(t),\int_t^T\tilde{\nu}(ds)\big{\ran} dt,
 \end{array}
$$
(where $\tilde{\nu}(\{0\})=0$, $\tilde{\nu}(\{T\})=0$ and $\tilde{\nu}(\Delta)=\nu(\Delta)$
for $\Delta\subset (0,T)$), and
$$
\int_0^T\lan h(t),q(t)\ran dt= \big{\lan} h(0),\int_0^T q(t)\big{\ran} dt +\int_0^T\big{\lan}\dot h(t),\int_t^Tq(s)\big{\ran} ds.
$$
Thus we can rewrite (2) as follows
$$
\begin{array}{l}
\big{\lan} h(0),\nu(\{0\})+\dis\int_0^T\tilde{\nu}(dt)\big{\ran} +  \lan h(T),\nu(\{T\})\ran-\big{\lan} h(0),\dis\int_0^T q(t)\big{\ran} dt\\
\qquad\qquad\qquad\qquad\qquad+ \dis\int_0^T\big{\lan}\dot h(t),p(t)+\dis\int_t^T\tilde{\nu}(ds)-\int_t^Tq(s)ds\big{\ran} dt=0
\end{array}
\eqno (4)
$$
Applying the equality for $h(\cdot)$ equal zero at the ends of the interval, we deduce that
$$
p(t)+\dis\int_t^T\tilde{\nu}(ds)-\int_t^Tq(s)ds= {\rm const} = c\quad  a.e.
$$
Changing $p(\cdot)$ on a set of  measure zero, if necessary, we get from here that $p(\cdot)$ is a function of bounded variation continuous from the left.  If we now apply (4) for $h(\cdot)$ equal to zero at zero, we conclude that $c= - \nu(\{T\}) $, that is
$$
p(t)+\dis\int_t^T\nu(ds)-\int_t^Tq(s)ds= 0,\quad\forall \; t>0.
\eqno (5)
$$
and (4) eventually implies that
$$
\int_0^T\nu(dt)=\int_0^Tq(t)dt.
\eqno (6)
$$

Defining $p(0)$ by continuity, we can conclude with the following statement..

\begin{proposition}\label{pro2}
If $(\xb(0),0,\ldots,0)$ is a local minimum of $J$ on $W^{1,1}\times \R_+^k$, then there are a measure
$\nu\in\sd\vf(\xb(\cdot))$, a function $p(\cdot)$ of bounded variation continuous from the left 
and a summable $q(\cdot)$ taking values in $\sd_C\lan p(t),\psi_0(t,\cdot)\ran(\xb(t))$ a.e. such that
$p(0)= \nu(\{0\})$ and the relations (3), (5) and 
(6) hold true. 
\end{proposition}



\subsection{Statement of the maximum principle}
Let us return to ({\bf OC}). To make the reduction theorem applicable, we have to be sure that
the collection $\cu$ of ``good" controls is sufficiently rich. To this end we add the following
assumption to (H$_1$)-(H$_5$):

\vskip 1mm

(H$_6$)  for almost every $t$ the mapping $f(t,\cdot,u)$ is Lipschitz near $\xb(t)$ for every $u\in  U(t)$, that is there are $\rho>0$ and $\ep >0$ (depending on $t$ and $u$) such that
$$
\| f(t,x,u) -f(t,x',u)\|\le \rho\| x-x'\|,\; \forall\; x,x'\in B(\xb(t),\del).
$$ 
This is the weakest Lipschitz-type assumption of $f(t,\cdot,u)$, certainly sufficient for applications. Note however that there are proofs of the maximum principle with $f(t,\cdot,u)$
Lipschitz at $\xb(t)$ only for $u=\ub(t)$ (see e.g \cite{AV04}), although  at the expense of some other assumptions. We for instance do not need continuity of $f(t,x,\cdot)$. It is not clear whether these two weakenings can be combined.

To state the maximum principle for ({\bf OC}) we recall the  standard notation: for a vector $p\in\R^n$ let
$$
H(t,x,p,u)= \lan p,f(t,x,u)\ran,\quad \ch(t,x,p) = \sup_{u\in U(t)}H(t,x,p,u).
$$
We also set $\overline{\Delta} = \{t:\; g(\xb(t))=0\}$ and
$$
\sd_C^>g(t,x)= \conv\{\lim y_m:\ y_m\in\sd g(t_m,\cdot)(x_m),\; t_m\to t,\ x_m\to x,\; g(t_m,x_m)>g(t,x)\}.
$$

\begin{theorem}\label{th1} Assume (H$_1$)-(H$_6$). If $(\xb(\cdot),\ub(\cdot)$ is a strong local minimum in
	({\bf OC}), then there are $\la\in [0,1]$, a function of bounded variation $p(t)$ on $[0,T]$, continuous   from the left,  a regular nonnegative measure $\mu$ with $\mu([0,T])\le 1$ supported on $\overline{\Delta}$, a summable
	$q(t)\in\sd_CH(t,\cdot,p(t),\ub(t))(\xb(t))$ a.e.,
	a $\mu$-measurable selection $\ga(t)$ of the set-valued mapping $t\to\sd_C^>g(t,\xb(t))$ and a pair  $(w_0,w_T)\in \la\sd\ell(\xb(0),\xb(T)) + \Phi'^*(\xb(0),\xb(T))\big(N(S,\Phi(\xb(0),\xb(T)))\big)$
	such that the following relations are satisfied
	
	\vskip 1mm
	
	\quad $\la+ \| p(\cdot)\|+ \mu([0.T])>0$ \ (nontriviality);
	
	\vskip 1mm
	
	\quad $p(0)=w_0,\quad  p(T)+\ga(T)\mu(\{T\})=- w_T$ \ (transversality);
	
	\vskip 1mm
	
	\quad $p(t)= -w_T +\dis\int_t^T q(s)ds -\dis\int_t^T \ga(s)d\mu(s),\;\forall\; t$ \   (adjoint equation); 
	
	\vskip 1mm
	
	\quad $H(t,\xb(t),p(t),\ub(t)) =  \ch(t,\xb(t),p(t))$ a.e. \  (maximim principle).
	
\end{theorem}
\noindent  Here $\| p(\cdot)\|=\sup_t\|p(t)\|$.



\begin{remark}\label{rem0}{\rm  In the classical smooth setting the statement of the theorem  reduces
to the maximum principle for optimal control problems proved in \cite{IT}. It differs however from the statement of the maximum principle in \cite{RV}, Theorem 9.3.1 (proved  under basically the same assumptions as here). Nonetheless,  both statements are equivalent. To see this, let us denote by $q_1(\cdot)$ what is $p(\cdot)$ in our theorem and by $\eta(t)$ what is $q(t)$, that is
$q_1(t)= -w_T +\int_t^T \eta (s)ds -\int_t^T \ga(s)d\mu(s)$,  set further $p(t)= w_0-\int_0^t\eta(s)ds$ and  $q(t) = p(t)+\int_0^t\ga(s)\mu(ds)$.  Then
$$
q_1(t)-q(t) = -(w_0+w_T)+\int_0^T\eta(T)dt -\int_0^T\ga(T)\mu(dt)=0        
$$
for all $t<T$.  Since $q_1$ is continuous from the left, we get from here that
$\lim_{t\to T}q(t) = q_1(T)$ and therefore $q(T)= q_1(T)+\ga(T)\mu(\{T\})= -w_T$.
Verification that Theorem 3.3 and Theorem 9.3.1 in \cite{RV} are equivalent is now straightforward .}
\end{remark}

\vskip 1mm

\subsection{Proof of the maximum principle}

1. The key element of the proof is the study of necessary conditions in ({\bf OC}$_k$) with
$U_{0}(t)\equiv \{\ub(t)\}$, so that $U_k(t)=\{\ub(t),u_1(t),\ldots,u_k(t)\}$. 
We claim that the theorem is valid if the following weaker integral maximum principle holds for every such ({\bf OC}$_k$). 

\begin{theorem}\label{th2} Assume (H$_1$)-(H$_6$). If $(\xb(\cdot),\ub(\cdot))$ is a strong local minimum in ({\bf OC$_k$}) with $U_k(t)=\{\ub(t),u_1(t),\ldots,u_k(t)\}$, $u_i(\cdot)\in\cu$,  then the conclusion of Theorem \ref{th1}
holds with the maximum principle replaced by
$$
\int_0^TH(t,\xb(t),p(t),\ub(t))dt \ge \dis\int_0^TH(t,\xb(t),p(t),u_i(t))dt,\; i=1,\ldots,k.
$$	
\end{theorem}

So assume that Theorem \ref{th2} is true. Given an ({\bf OC}$_k$), let $\Lambda_k$ be the set of triples $(\la,p(\cdot),\mu)$ satisfying the  conditions of Theorem \ref{th2}  along with normalized  nontriviality  condition: $\la+\|p(\cdot)\|+\mu([0,T])=1$. As all four relations in the theorem are positively homogeneous with respect to $\la,p(\cdot),\mu$, this set is nonempty by the  theorem.
It is clear that $\Lambda_{k'}\subset \Lambda_k$ if $U_k\subset U_{k'}$. 
If the intersection $\Lambda$ of all $\Lambda_k$ is nonempty, then for any $(\la,p(\cdot),\mu)\in\Lambda$
$$
\int_0^TH(t,\xb(t),p(t),\ub(t))dt \ge \dis\int_0^TH(t,\xb(t),p(t),u(t))dt,\; \forall\;
u(\cdot)\in \cu.
\eqno (7)
$$
In view of Proposition \ref{mea}, to prove the maximum principle we only need to
verify that (7) actually holds for all measurable selections of $U(\cdot)$ for which
the integral on the write side of (7) makes sense. Let $v(\cdot)$ be such a selection
of $U(\cdot)$. If the inequality opposite to (7) holds for this $v(\cdot)$, then
the set $\Delta=\{t:\; H(t,\xb(t),p(t),v(t))>H(t,\xb(t),p(t),\ub(t))\}$ has positive measure.  
Consider on $\Delta$ the function
$$
\eta_{\ep}(t) = \sup_{x,x'\in B(\xb(t),\ep),\ x\neq x'}\frac{\|f(t,x,p(t),v(t))-f(t,x',p(t),v(t))\|}{\|x-x'\|},
$$
which is measurable by (H$_3$). By (H$_6$) $\lim_{\ep\to 0}\eta_{\del}(t)<\infty$
for almost every $t\in \Delta$. It follows that there is a positive $\ep>0$ and a
subset $\Delta'\subset \Delta$ of positive measure on which $\eta_{\ep}$ is summable.
Then the $u(\cdot)$ coinciding with $v(t)$ on $\Delta'$ and with
$\ub(t)$ on the complement of $\Delta'$ obviously belongs to $\cu$ and (7) fails for this $u(\cdot)$.

Thus to prove the claim we have to verify that $\Lambda\neq\emptyset$. This in turn will
be true if we verify  that any $\La_k$ is compact if  measures are considered in the weak$^*$-topology.  
So fix a $\Lambda_k$, take  $(\la_m,p_m(\cdot),\mu_m)\in\Lambda_k$, $m=1,2,\ldots$, and let  $q_m(\cdot)$ and  $\ga_m(\cdot)$
be measurable selections of  $\sd_C\lan p_m(t),f(t,\cdot,\ub(t))\ran(\xb(t))$ and     $\sd_C^>g(t,\cdot)(\xb(t))$ respectively such that $p_m(t)= - w_{Tm}+\int_t^Tq_m(s)ds-\int_t^T\ga_m(s)d\mu_m(s)$ and $(p_m(0),-w_{Tm})\in\sd\ell(\xb(0),\xb(T) +A(N(S,\zb))$, where we have set for simplicity
$\zb=\Phi(\xb(0),\xb(T))$ and $A=\Phi'^*(\xb(0),\xb(T))$.
We may assume that $\la_k$ converge to some $\la\ge 0$ and $p_m(0)$ converge to some $p_0$.
Since the set of nonnegative measures with $\mu([0,T])\le 1$ is compact in the weak$^*$ topology of measures, we may also assume  (taking if necessary a subsequence) that
$\mu_k$ weak$^*$ converge
to some $\mu\ge 0$. In particular,  $\mu_m([0,T])\to \mu([0,T])$.

By (H$_2$), (H$_4$) $\| q_m(t)\|\le k_0(t)\| p(t)\|$, and $\|\ga_m(t)\|\le \rho_g$ almost everywhere.
It follows that $(w_{Tm})$ is a bounded sequence, and again we can assume that it converges to some $w_T$. 
Clearly $(p_0,w_T)\in \sd\ell(\xb(0),\xb(T))+A(N(S,\zb))$.  It also follows that the sequence of $q_m(\cdot)$ is weak compact in $L^1$
so that we may assume that it weakly converges to some $q(\cdot)$ and, consequently, the
integrals $\int_t^T q_m(s)ds$ converge uniformly to $\int_t^Tq(s)ds$. 

Finally, as
$\sd_C^>g(t,\cdot)(\xb(t))$ are convex and closed, there is a measurable $\ga(t)$ with values in
$\sd_C^>g(t,\cdot)(\xb(t))$ such that the measures $\nu$ with $d\nu_m(t)=\ga_m(t)d\mu_m(t)$
weak$^*$ converge to $\nu$ defined by $d\nu(t)=\ga(t)d\mu(t)$ (see \cite{RV}, Proposition 9.2.1). It follows that $p_m(t)$ converge to $p(t)=-w_T +\int_t^Tq(s)ds-\int_t^T\ga(s)\mu(ds)$
 at every $t>0$ at which $\int_t^T\ga(s)d\mu(s)$ is continuous.
 Setting $p(0)=p_0$ we find that the limiting objects satisfy all required relations except maybe nontriviality. To prove nontriviality of the limiting objects,
 observe that  if $\mu=0$, then $p_m(\cdot)$ converge
to $p(\cdot)$ uniformly and hence $\la + \| p(\cdot\|=1$. This completes the proof of the claim.   

\vskip 1mm

2. So we have to prove Theorem \ref{th2}. We shall do this by applying Theorem \ref{red2}
with $U_0(t)=\{\ub(t)\}$. It follows from the theorem 
(if we replace $\la$ by $K\la$ in case when $K>1$)
that there is a sequence
$(\xb_m(\cdot),\ub_m(\cdot))\in X_k$ with $\xb_m(\cdot)$ converging to $\xb(\cdot)$ uniformly
and $\ub_m(\cdot)$ converging to $\ub(\cdot)$ in $L^1$ such that  $(\xb_m(\cdot),0,\ldots,0)$ is a local minimum of the functional
$$
\begin{array}{l}
\la\psi_m(x(\cdot)+ d(\Phi(x(0),x(T)),S)
+J_{k}(x(\cdot),\ub_m(\cdot),\al_1,\ldots,\al_k)\\
\qquad\qquad\qquad
+ \  m^{-1}\Big(\|x(\cdot) -\xb_m(\cdot)\|_C +\dis\int_0^T\big( \sum\al_i\| u_i(t)-\ub_m(t)\|\big)dt         \Big)
\end{array}
$$
on $C([0,T])\times \R_+^k$ for any $m$. The structure of the functional  precisely corresponds to the basic model. To see this, it is enough to set
$$
\begin{array}{l}
\vf(x(\cdot))= \la\psi_m(x(\cdot)+ d(\Phi(x(0),x(T)),S)+m^{-1}\| x(\cdot)-\xb_m(\cdot)\|_C\\
\psi_0(t,x)= f(t,x,u_m(t));\\
\psi_i(t,x) = f(t,x,u_i(t))-f(t,x,\ub_m(t)) + m^{-1} \| u_i(t)-\ub_m(t)\|
\end{array}
$$
Thus we can apply Proposition \ref{pro2} to get a necessary condition
for  $(\xb_m(\cdot),0,\ldots,0)$ to be a local minimum of the functional. All we need is to find  suitable expressions  for elements of the subdifferential of $\vf$ using
the standard rules of subdifferential calculus collected in  the second section.

If $\nu\in\sd\psi_m(x(\cdot))$, then $\nu=\xi_1\nu_1+\xi_2\nu_2$, where $\xi_i\ge 0$,
$\xi_1+\xi_2=1$, 
$$
\xi_1(\ell(x(0),x(T))+m^{-2}-\psi_m(x(\cdot)))=\xi_2(\max_{t}g(t,x(t))-\psi_m(x(\cdot)))=0
$$
and $\nu_1$ is supported on $\{0,T\}$ with weights belonging to $\sd\ell(x(0),x(T))$, 
while $\nu_2(dt)=\ga(t)\mu(dt)$, where $\mu$ is a probability measure supported on the set
of $t$ at which $g(t,x(t))$ attains its maximum and $\ga(t)$ is a measurable selection of the set-valued mapping $t\to \overline{\sd}g(t,\cdot)(x(t))$. 

Setting $\la_1=\la\xi_1$ and writing $\mu$ instead of $\la\xi_2\mu$, we conclude that
for any $m$ we can find 
$\la_{m1}\ge 0$, a nonnegative measure $\mu_m$ supported on the set
$\Delta(\xb_m(\cdot))=\{t: g(t,\xb_m(t))=\psi_m(\xb_m(\cdot))\}$
 a measurable selection $\ga_m(\cdot)$ of the set-valued mapping $\overline{\sd}g(t,\cdot)(\xb_m(t))$,   a pair
$(w_{m0},w_{mT})$ of vectors in $\R^n$, a measurable selection $q_m(\cdot)$
of the set-valued mapping $t\to \sd_Cf(t,\cdot.\ub_m(t))(\xb_m(t))$ 
and a function $p_m(t)$ of bounded variation 
such that

\vskip 1mm

$\bullet$ $\la_{m1}+ \|\mu_m\|=\la$, \ $\la_{m1}(\ell(\xb_m(0),\xb_m(T))+m^{-2}-\psi_m(\xb_m(\cdot)))=0$ and \\ 
$\mu=0$ if $\max_tg(t,\xb_m(t))<\psi_m(\xb_m(\cdot) )$;

\vskip 1mm

$\bullet$ $(w_{m0},w_{mT})\in \la_{m1}\sd\ell(\xb_m(0),\xb_m(T))$

\vskip 1mm
\qquad\qquad\qquad\qquad\qquad\qquad $+\Phi'^*(\xb_m(0),\xb_m(T))\big(\sd d(\cdot,S)(\Phi(\xb_m(0),\xb_m(T)))\big)$

\vskip 1mm

\noindent and the following three relations are satisfied up to terms of order $m^{-1}$:
 
$$
\begin{array}{l}
p_m(t)+w_{mT}+\dis\int_t^T(\ga_m(s)\mu_m(ds)-q_m(s)ds)\ = \ 0\; {\rm a.e.};\\
w_{m0} + w_{mT} + \ \dis\int_0^T(\ga_m(t)\mu_m(dt)-q_m(t)dt)\ = \ 0,\\
\dis\int_0^T\lan p_m(t),\psi_i(t,\xb_m(t))\ran dt \ \le 0,\; i=1,\ldots,k.
\end{array}
$$
It follows that
 $p_m(0)=w_{m0} +O(m^{-1})$ and $p_m(T)= -(w_{mT}+\ga_m(T)\mu_m(\{T\}))+ O(m^{-1})$. Taking into account that
$\lan p,\psi_0(t,x)\ran  =H(t,x,p,\ub_m(t))$
and $\lan p,\psi_i(t,x)\ran =H(t,x,p,u_i(t))-H(t,x,p,\ub_m(t))$, we conclude
that $p_m(\cdot)$ satisfies the following system of relations for $i=1,\ldots,k$
valid up to terms of order $m^{-1}$:
$$
\begin{array}{l}
p_m(t)=p_m(T)+\ga_m(T)\mu(\{T\}+\dis\int_t^T\big(q_m(t)dt-\ga_m(s)\mu_m(ds)\big);\\ 
\dis(p_m(0),-p_m(T)+\ga_m(T)\mu(\{T\})\in \sd\la_{m1}\ell(\xb_m(0),\xb_m(T)) \\
\qquad\qquad\qquad\qquad\qquad +\Phi'^*(\xb_m(0),\xb_m(T))\big(\sd d(\cdot,S)(\Phi(\xb_m(0),\xb_m(T)))\big);\\
\dis\int_0^T\big(H(t,\xb_m(t),p_m(t),u_i(t))-H(t,\xb_m(t),p_m(t),\ub_m(t))   \big)dt\le 0.
\end{array}
\eqno (11)
$$

3. We can now easily finish the proof.   
By (H$_2$) the functions $\|q_m(t)\|$ are bounded by the same summable function, hence
uniformly integrable, by 
(H$_4$) $\|\ga_m(t)\|\le \rho_g$ a.e.,  the sequences of $p_m(0)$ and $w_{mT}=-(p_m(T)+\ga_m(T)\mu_m(\{T\}))$ are bounded, hence the sequence of $p_m(T)$ is also bounded. 
We may assume that $q_m(\cdot)$ weak converge to some $q(\cdot)$, hence
$\int_t^Tq_m(s)ds$ uniformly converge to $\int_t^Tq(s)ds$, and (as we obviously can assume that $\la_{m1}$ converge to some $\la$)   the pairs $(p_m(0),w_{mT})$ 
 converge to some  $(w_0,w_T)\in\la\sd\ell(\xb(0),\xb(T))+\sd d(\cdot,S)(\xb(0),\xb(T))$. 

The sequence of measures $ \ga_m(t)\mu_m(dt)$ is weak$^*$ compact.
We refer again to  Proposition 9.2.1 of \cite{RV} to conclude that
the limiting measures have the form $\ga(t)\mu(dt)$, where $\mu$ is a weak$^*$ limit of $\mu_m$
and $\ga(t)\in \cap_m\cl \conv(\cup_{s>m}\overline{\sd}_Cg(t,\cdot)(\xb_s(t))$ almost everywhere.
As a result, we deduce  (as $\| p_m(\cdot)\|$ are uniformly bounded) that $p_m(\cdot)$
converge almost everywhere to $p(t)= - w_T  + \int_t^T(q(s)-\ga(s)d\mu(s))$ and $p(0)=w_0$, $p(T)+\ga(T)\mu({T}) = -w_T$. The transversality condition,  the adjoint inclusion
and the integral inequality of Theorem \ref{th2} is now immediate from (11). The nontriviality
condition is obvious in the non-singular case as $\la>0$. To prove that this condition
holds also in the singular case, with $\la =0$, we have to recall that $\xb_m(\cdot)$ are not feasible
in ({\bf OC}$_k$) which means that $\yb_m=\Phi(\xb_m(0),\xb_m(T))\not\in S$. Therefore the norm of any
element of $\sd d(\cdot,S)(\yb_m)$ is $1$. Taking into account that $\mu=0$
and $\Phi'^*(\xb_m(0),\xb_m(T))$ is one-to-one by (H$_5$), we conclude that $\|(p(0),-p(T))\|>0$.

 Finally, it  is an easy matter to see that  the set
 $\cap_m\cl\conv(\cup_{s>m}\overline{\sd}_Cg(t,\cdot)(\xb_s(t))$ 
 lies in $\sd_C^>g(t,\cdot)(\xb(t))$.
Indeed, it is sufficient to note that (in the non-singular case) $J_m(\xb_m(\cdot))>0$, and therefore  $\max_t g(t,\xb_m(t))>0$.
(Otherwise $\xb_m(\cdot)$ would be admissible in ({\bf OC}) and $\ell(\xb_m(0),\xb_m(T))<0=\ell(\xb(0),\xb(T))$.)
This completes the proof of the theorem.

\section{Second order conditions}
Below we use the following notation. Let $f$ be a function or a mapping depending on two vector variables, e.g. $f(x,u)$. Then $f'$ stands for the full derivative of $f$, while the partial derivatives will be denoted $f_x'$ and $f_u'$, the same with second derivatives.  

\subsection{Basic unconstrained model}
As in the first order theory we start here with introduction of a certain basic unconstrained model. This time this is the problem of minimizing
$$
 J(x,u,\al_1,\ldots,\al_k)=g(F_0(x,u)+\sum_{i=1}^k\al_iF_i(x,u))
\eqno (12)
$$
subject to $x\in X$, $u\in U,\; \al_i\ge 0$. Here $X$ is a Banach space, $U$ is a closed subset in a Banach space $W$ and $F_i$ are mappings from $X\times W$ into a Banach space $Y$.
As usual we fix some $(\xb,\ub)\in X\times U$ and assume that $J$ has a local minimum on 
$X\times U\times \R_+^k$ at $(\xb,\ub,0,\ldots,0)$. We further assume that

\vskip 1mm

(H$_7$)  
$g$ is  continuous and sublinear;
the mappings $F_i$, $i=0,\ldots,k$,  are continuous  and continuously differentiable  near $(\xb,\ub)$ and  twice differentiable at $(\xb,\ub)$.        

\vskip 1mm

Set $\yb=F_0(\xb,\ub)$. 
Since the function $(x,\al_1,\ldots,\al_k) \to J(x,\ub,\al_1,\ldots,\al_k)$ attains a local minimum on $X\times \R_+^k$ at $(\xb,0,\ldots,0)$, there is a $y^*\in \sd g(\yb)$ such that
$$
(y^*\circ F_0)_x'(\xb,\ub)=0,\quad (y^*\circ F_i)(\xb,\ub)\ge 0
\eqno (13)
$$
We shall denote by $\Lambda$ the collection of $y^*\in\sd g(\yb)$ satisfying (13).

 Set $\mbox{\boldmath $x$}=(x,\al_1,\ldots,\al_k)$,
$\mbox{\boldmath $\xb$}=(\xb,0,\ldots,0)$,
and let $F(\mbox{\boldmath $x$},u)$ stand for $F_0(x,u)+\sum\al_iF_i(x,u)$. 
We 
refer to $\mbox{\boldmath $x$}$ as {\it feasible} if $\al_i\ge 0$.
Note first that for any $\ep>0$ there is an $\eta>0$ such that
$$
g(F(\mbox{\boldmath $x$},u)+ F_{\mbox{\boldmath $x$}}'(\mbox{\boldmath $\xb$},\ub)\mbox{\boldmath $h$})\ge g(F(\mbox{\boldmath $x$}+\mbox{\boldmath $h$},u)- (\ep/2)\| \mbox{\boldmath $h$}\|\ge g(\yb)- (\ep/2)\|\mbox{\boldmath $h$}\|,
\eqno(14)
$$
if $\|\mbox{\boldmath $x$}- \mbox{\boldmath $\xb$}\|\le\eta$, $\|\mbox{\boldmath $h$}\|\le\eta$,
$\| u-\ub\|\le\eta$. Here we have set $\mbox{\boldmath $h$}=(h,\beta_1,\ldots,\beta_k)$,
$\| \mbox{\boldmath $x$}\|=\| x\|+\sum|\al_i|$, $\| \mbox{\boldmath $h$}\|=\| h\|+\sum|\beta_i|$. The left inequality is valid for all such $\mbox{\boldmath $x$},\mbox{\boldmath $h$}\in X\times \R_k$ while the right one, of course, only if $\mbox{\boldmath $x$},\mbox{\boldmath $h$}\in X\times\R_+^k$.

Indeed, fix some $\mbox{\boldmath $x$},\ u,\ \mbox{\boldmath $h$}$ and choose a $y^*\in\sd g(0)$ such that $g(F(\mbox{\boldmath $x$}+\mbox{\boldmath $h$},u))=\lan y^*,F(\mbox{\boldmath $x$}+\mbox{\boldmath $h$},u)\ran$. 
Then
$$
g(F(\mbox{\boldmath $x$},u)+ F_{\mbox{\boldmath $x$}}'(\mbox{\boldmath $\xb$},\ub)\mbox{\boldmath $h$})- g(F(\mbox{\boldmath $x$}+\mbox{\boldmath $h$}),u)
\ge\lan y^*,F(\mbox{\boldmath $x$},u)+F_{\mbox{\boldmath $x$}}'(\mbox{\boldmath $\xb$},\ub)\mbox{\boldmath $h$})- F(\mbox{\boldmath $x$}+\mbox{\boldmath $h$},u))\ran.
$$
But $\| F(\mbox{\boldmath $x$},u)+F_{\mbox{\boldmath $x$}}'(\mbox{\boldmath $\xb$},\ub)\mbox{\boldmath $h$}- F(\mbox{\boldmath $x$}+\mbox{\boldmath $h$},u))\| =r(x,u,\mbox{\boldmath $h$})\| \mbox{\boldmath $h$}\|$, 
where $r(\mbox{\boldmath $x$},u,\mbox{\boldmath $h$})\to 0$ when $\mbox{\boldmath $x$}\to\mbox{\boldmath $\xb$}$, $u\to\ub$,
$ \mbox{\boldmath $h$}\to 0$ and  $\| y^*\|$ does not exceed the Lipschitz constant of $g$.

It will be convenient to set for  further discussions  
$$
p_{\ep}(\mbox{\boldmath $x$},u,\mbox{\boldmath $h$})= 
g(F(\mbox{\boldmath $x$},u)+ F_{\mbox{\boldmath $x$}}'(\mbox{\boldmath $\xb$},\ub)\mbox{\boldmath $h$})
+\ep\| \mbox{\boldmath $h$}\|.
$$ 
By (14) $p_{\ep}(\mbox{\boldmath $x$},u,\mbox{\boldmath $h$})\ge g(F(\mbox{\boldmath $\xb$},\ub))$ for feasible  $\mbox{\boldmath $x$},\mbox{\boldmath $h$}$
satisfying $\|\mbox{\boldmath $x$}-\mbox{\boldmath $\xb$}\|\le\eta$,
$\|\mbox{\boldmath $h$}\|\le\eta$, $\| u-\ub\|\le\eta$.
%
We claim that for every $\ep>0$ there are $\eta>0$ and $\del>0$
such that
$$
\inf_{\mbox{\boldmath $h$}\in X\times \R_+^k} 
p_{\ep}(\mbox{\boldmath $x$},u,\mbox{\boldmath $h$})
=\inf_{\mbox{\boldmath $h$}\in X\times\R_+^k,\ \| \mbox{\boldmath $h$}\|\le\dis\eta}
p_{\ep}(\mbox{\boldmath $x$},u,\mbox{\boldmath $h$})\ge g(F(\mbox{\boldmath $\xb$},\ub))
\eqno (15)
$$
whenever $\|\mbox{\boldmath $x$}-\mbox{\boldmath $\xb$} \|\le \del$, $\|u-\ub\|\le\del$.  
 Indeed,
taking a smaller $\eta$, if necessary, we can be sure that $g(F(\mbox{\boldmath $x$}+\mbox{\boldmath $h$},u))\ge g(F(\mbox{\boldmath $\xb$},\ub))$
if $\| \mbox{\boldmath $x$}-\mbox{\boldmath $\xb$}\|<\eta$, $\| u-\ub\|<\eta$,  $\|\mbox{\boldmath $h$} \|<\eta$.  
Then by (14) $p_{\ep}(\mbox{\boldmath $\xb$},\ub,0)=g(\yb)$ and $p_{\ep}(\mbox{\boldmath $\xb$},\ub,\mbox{\boldmath $h$})\ge g(\yb)+(\ep\eta)/2$ if $\|\mbox{\boldmath $h$}\|=\eta$. It remains to choose
$\del>0$ such that $|g(F(\mbox{\boldmath $x$},u)+ F_{\mbox{\boldmath $x$}}'(\mbox{\boldmath $\xb$},\ub)\mbox{\boldmath $h$})-g(\yb)|<\ep\eta/4$ if $\|\mbox{\boldmath $x$}-\mbox{\boldmath $\xb$}\|<\del$, $\| u-\ub\|<\del$  and $\|\mbox{\boldmath $h$}\|\le \eta$.  As $p_{\ep}(\mbox{\boldmath $x$},u,\cdot)$ is a convex function, it follows
that its lower bound is realized in the $\eta$-ball around zero and (15) follows.

Let $\Lambda_{\ep}$ and $\Lambda_{0\ep}$ be the sets  of 
$y^*\in\sd_{\ep} g(\yb)$ and $y^*\in\sd g(0)$ respectively
satisfying 
$$
\| (y^*\circ F_0)_x'(\xb,\ub))\|\le \ep,\quad \lan y^*,F_i(\xb,\ub)\ran \ge  -\ep,\; i=1,\ldots,k.
$$ 
(Here $\sd_{\ep}$ stands for the $\ep$-subdifferential in the sense of convex analysis:
$y^*\in\sd_{\ep}g(y)$ if $g(z)-g(y\ge \lan y^*,z-y\ran -\ep$ for all $z$.)
Our final claim is that we can choose  $\del>0$  such that in addition to (15)
	$$
	\inf_{\mbox{\boldmath $h$}\in X\times \R_+^k} p_{\del}(\mbox{\boldmath $x$},u,\mbox{\boldmath $h$})=\sup_{y^*\in\La_{\ep}}\lan y^*,F(\mbox{\boldmath $x$},u)).
	\eqno (16)
	$$
Indeed, 
$$
\begin{array}{lcl}
\dis\inf_{\mbox{\boldmath $h$}}p_{\ep}(\mbox{\boldmath $x$},u,\mbox{\boldmath $h$})&=&\dis\inf_{\mbox{\boldmath $h$}}
(\dis\sup_{y^*\in\sd g(0)}\lan y^*, F(\mbox{\boldmath $x$},u)+F_{\mbox{\boldmath $x$}}'(\mbox{\boldmath $\xb$},\ub)\mbox{\boldmath $h$}\ran + \ep\|\mbox{\boldmath $h$}\|)\\
&= &\dis\sup_{y^*\in\sd g(0)}(\lan y^*, F(\mbox{\boldmath $x$},u)\ran         +\dis\inf_{\mbox{\boldmath $h$}}(\lan y^*,F_{\mbox{\boldmath $x$}}'
(\mbox{\boldmath $x$},\ub)\mbox{\boldmath $h$}\ran 
+ \ep\|\mbox{\boldmath $h$}\|  ))   )\\
&=&\dis\sup_{y^*\in\Lambda_{0\ep}}(\lan y^*, F(\mbox{\boldmath $x$},u)\ran +\dis\inf_{ \mbox{\boldmath $h$}}(\lan y^*,F_{\mbox{\boldmath $x$}}'(\mbox{\boldmath $\xb$},\ub) \mbox{\boldmath $h$}\ran + \ep\|\mbox{\boldmath $h$}\|  )) \\
&=& \dis\sup_{y^*\in\Lambda_{0\ep}}\lan y^*, F(\mbox{\boldmath $x$},u)\ran
= \dis\sup_{y^*\in\Lambda_{\ep}}\lan y^*, F(\mbox{\boldmath $x$},u)\ran.
\end{array}
$$
The first equality follows from definitions. The second equality follows from the standard
minimax theorem (thanks to the fact that $\sd g(0)$ is weak$^*$-compact). To justify the third  equality we recall that
$$
\lan y^*,F_{\mbox{\boldmath $x$}}'(\mbox{\boldmath $\xb$},\ub)\mbox{\boldmath $h$}\ran=  \lan(y^*\circ F_0)_x'(\xb,\ub),h\ran + \lan y^*,\sum\beta_iF_i(\xb,\ub)\ran
$$
so that the infimum over $\mbox{\boldmath $h$}=(h,\beta_1,\ldots,\beta_k)\in X\times\R_+^k$ is $-\infty$ if either\\ $\| (y^*\circ F_0)_x'(\xb,\ub)\|>\ep$ 
or $\lan y^*,F_i(\xb,\ub)\ran <-\ep $ for some $i$. Finally, if $y^*\in\sd g(0)$, does not belong to $\sd_{\ep} g(\yb)$, then $\lan y^*,\yb\ran\le g(\yb)-\ep$ 
(recall that $g$ is a sublinear function) and the last equality
follows from obvious chain of inequalities below  (where  $z^*\in\Lambda$ and $L$ is the Lipschitz constant of $g$)
$$
\begin{array}{l}
\lan z^*,F(\mbox{\boldmath $x$},u)\ran\ge\lan z^*,\yb\ran-L\|F(\mbox{\boldmath $x$},u)-\yb\|\\ \qquad\qquad\quad\;\;=g(\yb)-L\|F(\mbox{\boldmath $x$},u)-\yb\| 
\ge\lan y^*,F(\mbox{\boldmath $x$},u)\ran+\ep - 2L\| F(\mbox{\boldmath $x$},u )-\yb\|
\end{array}
$$
if $\del$ is so small that $\ep - 2L\| F(\mbox{\boldmath $x$},u )-\yb\|>0$
if $\|\mbox{\boldmath $x$}-\mbox{\boldmath $\xb$}\|<\del$,\ $\|u-\ub\|<\del$.

Combining (14)-(16), we conclude with

\begin{proposition}
Under the assumptions, for any $\ep >0$
$$
\begin{array}{l}
g(F_0(\xb,\ub))=g(F(\mbox{\boldmath $\xb$},\ub))\le 
\dis\sup_{y^*\in\Lambda_{\ep}}\lan y^*, F(\mbox{\boldmath $x$},u)\ran \\
\qquad\qquad\qquad\qquad\qquad\quad
=\dis\sup_ {y^*\in\Lambda_{\ep}}\lan y^*,F_0(x,u)+ \sum\al_iF_i(x,u)    \ran
\end{array}
$$
for all $(\mbox{\boldmath $x$},u)=(x,u,\al_1,\ldots,\al_k)\in X\times U\times \R_+^k$ close to 
$(\mbox{\boldmath$\xb$},\ub)=(\xb,\ub,0,\ldots,0)$.

\end{proposition}

We are ready now to state and proof the main result of this subsection. It will be done under the following additional assumption:

(H$_8$) $W$ is densely embedded into another Banach space $V$ and
$F_{0u}'(\xb,\ub)$ extends by continuity to the whole of $V$. 

Recall that by definition $W$ is densely embedded into $V$ if there is a one-to-one linear
mapping $i: W\to V$ such that $i(W)$ is dense in $V$ and $\| u\|_W\ge \|i(u)\|_V$. As usual we identify $W$ and $i(W)$.  








Finally, we need  the concept of a {\it critical set} of
$J$ at $(\xb,\ub)$ which is the collection of all tuples
$(h,u,\beta_1,\ldots,\beta_k)$ with $h\in X,\ u\in W,\ \beta_i\ge 0$  such that
$$
\begin{array}{l}
g\big(F_0(\xb,\ub)+F_0'(\xb,\ub)(h,u)+\dis\sum_{i=1}^k\beta_iF_i(\xb,\ub) \big)\le g(F_0(\xb,\ub)).
\end{array}
\eqno (17)
$$
We shall denote this set by $\crit J$

\begin{theorem}\label{th4}
Assume (H$_7$), (H$_8$).   Then

$$
\begin{array}{l}
\dis\sup_{y^*\in\Lambda}\lan y^*,F_0''(\xb,\ub)(h,u)(h,u)+2F_{0u}'(\xb,\ub)v\\
\qquad\qquad\qquad\qquad\qquad +\dis\sum_{i=1}^k\beta_i(F_{ix}'(\xb,\ub)h+ F_{iu}'(\xb,\ub)u)\ran \ge 0
\end{array}
$$
whenever $(h,u,\beta_1,\ldots,\beta_k)\in\crit J$,  $v\in T^2(U,\ub;u)$ (in $V$) and the following
property holds for $v$: for any $\ep>0$ there are $v_m\in W$ converging to $v$ in $V$ and $t_m\to 0$ such that $\ub+t_mu+t_m^2v_m\in U$ and 
$$
\begin{array}{l}
F_0(\xb+th,\ub+t_mu+t_m^2v_m)= F_0(\xb,\ub)+ tF_0'(\xb,\ub)(h,u) \\
\qquad\qquad\qquad\qquad\qquad+\dfrac{t_m^2}{2}(F_0''(\xb,\ub)(h,u)(h,u) +
2F_{0u}'(\xb,\ub)v_m) +  t_m^2r_m,
\end{array}
\eqno (18)
$$ 
where $\| r_m\|_V\le\ep$.
\end{theorem}
\noindent Here of course $(h,u)\to F_0''(h,u)(h,u)$ is the quadratic form associated with the
second order derivative of $F_0$ at $(\xb,\ub)$.

\begin{remark}\label{rem3}{\rm    \ As $g$ is a sublinear function, that is $g(\la y)=\la g(y)$
and $g(y+z)\le g(y)+g(z)$, the critical set contains the {\it critical cone}, obtained if we replace  (17) by
$$
\{(h,u,\beta_1,\ldots,\beta_k):\; g\big (F_0'(\xb,\ub)(h,u)+\dis\sum_{i=1}^k\beta_iF_i(\xb,\ub)\big )\le 0 \},
$$		
lies in ${\rm Crit} J$ and coincides with it if $F_0(\xb,\ub)=0$.}
\end{remark}

\proof  We may assume without loss of generality that $\| y^*\|\le 1$ for elements of $\sd g(0)$. 
 By Proposition 5.1 and (18), for sufficiently large $m$
$$
\begin{array}{lcl}
g(F(\mbox{\boldmath $\xb$},\ub))&\le&\dis\sup_{y^*\in\Lambda_{\ep}}
\lan y^*, F(\mbox{\boldmath $\xb$}+t_m\mbox{\boldmath $h$},\ub+t_mu+t_m^2v_m)\ran\\
&=&\dis\sup_{y^*\in \Lambda_{\ep}}\lan y^*, F(\mbox{\boldmath $\xb$},\ub)+
t_mF'(\mbox{\boldmath $\xb$},\ub)(\mbox{\boldmath $h$},u) \\ 
&&\qquad\quad+(t_m^2/2)(F''(\mbox{\boldmath $\xb$},\ub)(\mbox{\boldmath $h$},u)(\mbox{\boldmath $h$},u) +2F_{u}'(\mbox{\boldmath $\xb$},\ub)v_m + 2r_m\ran\\
& \le& g(F(\mbox{\boldmath $\xb$},\ub)+t_mF'(\mbox{\boldmath $\xb$},\ub)(\mbox{\boldmath $h$},u)) \\
&&\qquad\quad +(t_m^2/2)\dis\sup_{y^*\in\Lambda_{\ep}}\lan y^*,
F''(\mbox{\boldmath $\xb$},\ub)(\mbox{\boldmath $h$},u)(\mbox{\boldmath $h$},u) +2F_{u}'(\mbox{\boldmath $\xb$},\ub)v_m+ r_m\ran .
\end{array} 
$$

We note next that 
$g(F(\mbox{\boldmath $\xb$},\ub)+ t_mF'(\mbox{\boldmath $\xb$},\ub)(\mbox{\boldmath $h$},u))\le g(F(\mbox{\boldmath $\xb$},\ub))$ since
$(\mbox{\boldmath $h$},u)= (h,u,\beta_1,\ldots,\beta_k)\in {\rm Crit} J$.
Thus
$$
0\le \sup_{y^*\in\Lambda_{\ep}}\lan y^*,
F''(\mbox{\boldmath $\xb$},\ub)(\mbox{\boldmath $h$},u)(\mbox{\boldmath $h$},u) +2F_{u}'(\mbox{\boldmath $\xb$},\ub)v_m\ran + \ep.
$$
By (H$_8$) we can pass to limit as $m\to \infty$ and write $v$ instead of $v_m$.
The  inequality is valid for any $\ep>0$.  
To conclude the proof we  note that $\La_{\ep'}\subset \La_{\ep}$ if $\ep'<\ep$, the function under the sign
of supremum is weak$^*$-continuous with respect to $y^*$ and  every $\La_{\ep}$
is a weak$^*$-compact set (as a subset of $\sd g(0)$). Hence we can pass to the limit as $\ep\to 0$ and get
$$
0\le \sup_{y^*\in\Lambda}
\lan y^*,F''(\mbox{\boldmath $\xb$},\ub)(\mbox{\boldmath $h$},u)(\mbox{\boldmath $h$},u) +2F_{u}'(\mbox{\boldmath $\xb$},\ub)v\ran
$$
which is precisely what has been stated.  \endproof

\subsection{Back to optimal control} 

We shall consider the problem with more specialized end point constraints:
$$
\begin{array}{rl}
{\rm minimize}& \ell_0(x(0),x(T)),\\
{\rm s.t.}& \dot x= f(t,x,u),\quad u\in U(t);\\
&g_i(t,x(t))\le 0,\; i=1,\ldots,s;\\
& \ell_j(x(0),x(T))\le 0,\; j=1,\ldots,l,\\ & \ell_j(x(0),x(T))=0,\; j=l+1,\ldots,r.
\end{array}
\leqno \rm ({\bf OC2})
$$

Here $U(t)$ are closed subsets of an Euclidean space 
and $(\xb(\cdot),\ub(\cdot))\in W^{1,1}\times L^{\infty}$ is a strong local minimum in the problem.
For brevity we write $\fb (t)$ for $f(t,\xb(t),\ub(t))$ and, likewise,
$\fb'(t)$ and $\fb''(t)$ for derivatives of $f$ in $ (x,u)$ at $(\xb(t),\ub(t))$. 
 We   further assume   that 

\vskip 1mm

(H$_9$)  the functions $\ell_j$, $j=0,1,\ldots,r$ are continuous and continuously differentiable near $(\xb(0),\xb(T))$
and have second derivatives at $(\xb(0),\xb(T))$;

\vskip 1mm

(H$_{10}$)  \  $\fb(t)$, $\fb'(t)$ and $\fb''(t)$ are summable on $[0,t]$ and 
 there is a $\del>0$ such that  
for a.e. $t$
the mapping  $(x,u)\to f(t,x,u)$ is 

-- continuous on $B(\xb(t),\del)\times U(t)$  along with its derivative w.r.t. $x$. Moreover for any $K>0$
the functions $\sup\{\| f(t,\xb(t),u)\|:\; u\in U(t),\  \|u\|\le K  \}$ and $ \sup\{\| f_x'(t,\xb(t),u)\|:\; u\in U(t),\  \|u\|\le K  \}$ are summable;

-- continuously differentiable at points of $B(\xb(t),\del)\times (B(\ub(t),\del)\cap U(t))$;

--  twice differentiable at $(\xb(t),\ub(t))$ uniformly in $t$ in the sense that
$$
f(t,\xb(t)+\la h,\ub(t)+\la u)= \fb(t)+\la \fb'(t)(h,u)+
(\la^2/2)\fb''(t)(h,u),(h,u)+r_{\la}(t,h,u),
$$
with $\int_0^T\|r_{\la}(t,h,u)\|dt=o(\la^2)\| h\|\|u\|$.

(H$_{11}$) $g_i$ and their derivatives with respect to $x$ are continuous;  the second derivatives $g_{ix}''(t,\xb(t))$ exist and are continuous on $[0,T]$. 

\vskip 1mm


\begin{theorem}\label{th5}
	Assume (H$_9$)-(H$_{11}$). If $(\xb(\cdot),\ub(\cdot)$ is a strong local minimum in the problem, then there are numbers $\la_0,\ldots,\la_r$, a function of bounded variation $p(t)$ on
	$[0,T]$ and  regular nonnegative measures $\mu_i$ supported on the sets
	$\Delta_i=\{t:\; g_i(t,\xb(t))=0\}$ such that 
	$\la_i\ge 0$; $\la_i\ell_i(\xb(0),\xb(T))=0,\ i=1,\ldots l$, and the following conditions are satisfied:
	
	\qquad $\dis\sum_{0\le i\le l}\la_i+\dis\sum_{i=l+1}^r|\la_i|+ \dis\sum_{i=1}^s\mu_i([0,T])=1$; 
	
	\qquad $(p(0),-(p(T)+
	\dis\sum_{i=1}^sg_{ix}'(T,\xb(t)))\mu_i(\{T\}) = \dis\sum_{i=0}^r\la_i\ell_i'(\xb(0),\xb(T))$;

	\qquad $p(t)= - \Big(p(T)+
	\dis\sum_{i=1}^s g_{ix}'(T,\xb(t)))\mu_i(\{T\}\Big) + \dis\int_t^TH_x'(s,\xb(s),p(s),\ub(s))ds$
	
	\qquad\qquad\qquad\qquad\qquad\qquad \qquad\qquad\qquad\qquad\qquad $ - \dis\sum_{i=1}^m \dis\int_t^Tg_{ix}'(s,\xb(s))\mu_i(ds)$;

	\qquad  $H(t,\xb(t),p(t),\ub(t))= \ch(t,\xb(t),p(t))$, \ a.e..
\end{theorem}
\noindent We shall denote by $\Lambda$ the collection of all such
$(\la_0,\ldots,\la_r,\mu_1,\ldots,\mu _s,p(\cdot))$.

The theorem  is a consequence of  Theorem \ref{th1}. Indeed,   set 

\vskip 1mm

$S=\{\xi=(\xi_1,\ldots,\xi_r)\in \R^r, \xi_j\le 0,\ j=1,\ldots,l,\; \xi_j=0,\ j=l+1,\ldots,r\}$

\vskip 1mm

$\Phi(x,x')= (\ell_1(x,x'),\ldots, \ell_r(x,x'))$,

\vskip 1mm

\noindent and take into account the following three elementary observations. The first is that   either the derivatives $\ell_i'(\xb(0),\xb(T))$,
$i=l+1,\ldots,r$ are linearly dependent or there is a $K>0$ such that
$d(\Phi(x(0),x(T)),S)\le K\sum_{i=l+1}^r|\ell_i(x(0),x(T))|$ if $x(\cdot)$ is close to $\xb(\cdot)$.
The second observation is that by (H$_{10}$) for $g(t,x)=\max_ig_i(t,x)$ we have
$\sd_C^>g(t,\cdot)(x)= \conv\{g_{ix}'(t,x): \;i\in I(t,x)\}$, where $I(t,x)=\{i:\;  g_i(t,x)=g(t,x)\}$.
Finally, it is an easy matter to see that in case when all $\la_i=0$ and $\mu_i=0$, the only
$p(\cdot)$ that can satisfy the theorem is identical zero, so there is no need to include $p(\cdot)$ in the nontriviality condition.

We can now state the main result of this section. 
We shall assume in what follows that ({\bf OC2}) is non-singular
at $(\xb(\cdot),\ub(\cdot))$, specifically  that 

\vskip 1mm

(H$_{12}$) the equation 
$-\dot p = H_x'(t,\xb(t),p,\ub(t))$ does not have a solution satisfying  
$$
(p(0),-p(T))=\sum_{i=l+1}^k\la_i\ell_i'(\xb(0),\xb(T))\;\&\;  
H(t,\xb(t),p(t),\ub(t))= \ch(t,\xb(t),p(t))
\eqno (19)
$$
unless all $\la_i$, $i=l+1,\ldots,r$ are zeros and therefore $p(t)\equiv 0$.

Let $I=\{0\}\cup\{i\in\{1,\ldots,l\}: \; \ell_i(\xb(0),\xb(T))=0  \}$, and let $\Delta_i(\del)$ stand for the $\del$-neighborhood of $\Delta_i$. 
Given a bounded measurable selection $w(t)$ of $U(t)$,
we shall consider 
the {\it critical cone} $C(w(\cdot))$ {\it associated with} $w(\cdot)$  which is 
the collection of triples $(h(\cdot),u(\cdot),\beta)\in W^{1,1}\times L^{\infty}\times \R_+$
such that $u(\cdot)\in T(U(\cdot),\ub(\cdot))$ and 
$$
\begin{array}{l}
\ell_i'(\xb(0),\xb(T))(h(0),h(T))\le 0; \; i\in I;\\
\ell_i'(\xb(0),\xb(T))(h(0),h(T))= 0,\; i=l+1,\ldots,r;\\
g_i(t,\xb(t))+ g_{ix}'(t,\xb(t))h(t)\le 0,\;\forall\; t;\\ 
\dot h(t) = \fb_{x}'(t)h(t) + \fb_u'(t)u(t)\\
\qquad\qquad\qquad\quad\ 
+\beta(f(t,\xb(t),w(t))-f(t,\xb(t),\ub(t))).
\end{array}
\eqno (20)
$$
for some $\del>0$. 
 
Finally, we shall say
that $u(\cdot)\in L^{\infty}$ and a measurable $v(\cdot)$ form a pair of {\it second order feasible variations} if $u(t)\in T(U(t),\ub(t))$, $v(t)\in T^2(U(t,\ub(t);u(t)))$ a.e. and   
there are $\bar{\la}>0$ and a summable  $\xi(\cdot)$ such that for $0<\la\le \bar{\la}$
$$
\int_0^T\| \fb_u'(t)v(t)\|dt<\infty;\quad\| \fb_u'(t)\| d(\ub(t)+\la u(t),U(t)) \le\la^2 \xi(t)\quad a.e..
$$

\begin{theorem}\label{th6}
Assume (H$_9$)-(H$_{12}$). Then 
for any bounded measurable selection $w(\cdot)$ of $U(\cdot)$,  any $(h(\cdot),u(\cdot))\in C(w(\cdot))$ 
 and  any measurable $v(\cdot)$ such that $(u(\cdot),v(\cdot))$ form a pair of second order feasible variations   we can find  a collection of multipliers in $\Lambda$, say
  $(\la_0,\ldots,\la_r,\mu_1,\ldots,\mu_s,p(\cdot))$, such that
$$
\begin{array}{l}
\dis\sum_{i=1}^r\la_i\ell_i''(\xb(0),\xb(T))((h(0),h(T)),(h(0),h(T)))\\
\qquad\qquad\qquad\qquad\qquad\qquad+\dis\sum_{i=1}^s
\dis\int_0^Tg_{ix}''(t,\xb(t))(h(t),h(t))\mu_i(dt)\\
 \qquad
   -\dis\int_0^T\Big(H_{(x,u)}''(t,\xb(t),p(t),\ub(t))( (h(t),u(t))(h(t),u(t)))\\
\qquad\qquad\qquad \qquad\qquad\qquad\qquad\qquad
+2H_u'(t,\xb(t),p(t),\ub(t))v(t)\\  \\
\qquad\qquad\qquad + \beta\big((H_x'(t,\xb(t),p(t),\ub(t))-H_x'(t,\xb(t),p(t), w(t)))h(t)\\
 \qquad\qquad\qquad\qquad\qquad\qquad\qquad\quad
 +H_u'(t,\xb(t),p(t),\ub(t))u(t)\big)\Big)dt\ge 0.
 \end{array}
\eqno (21)
$$ 
\end{theorem}

\subsection{Some comments}
1. If $\beta =0$,  the inequality in (21) assumes the form
$$
\begin{array}{l}
\dis\sum_{i=1}^r\la_i\ell_i''(\xb(0),\xb(T))((h(0),h(T)),(h(0),h(T)))\\
\qquad\qquad\qquad\qquad\qquad\qquad+\dis\sum_{i=1}^s
\dis\int_0^Tg_{ix}''(t,\xb(t))(h(t),h(t))\mu_i(dt)\\
\qquad
-\dis\int_0^T\Big(H_{(x,u)}''(t,\xb(t),p(t),\ub(t))( (h(t),u(t))(h(t),u(t)))\\
\qquad\qquad\qquad \qquad\qquad\qquad\qquad\qquad
+2H_u'(t,\xb(t),p(t),\ub(t))v(t)\ge 0
\end{array}
\eqno (22)
$$
 This is an infinitesimal second order condition defined by
the behavior of the functions at points arbitrarily close to $(\xb(t),\ub(t))$.
So it may be  tempting to consider it a necessary condition for a weak minimum. 
But $\Lambda$
is the set of Lagrange multipliers for which the maximum principle holds.
So it would be interesting to find an example of  a weak minimum for which (22) does not hold.
The first relation (21), on the other hand, deals with controls that can be arbitrarily far from $\ub(\cdot)$.
This in a sense makes Theorem \ref{th6}  a ``real" second order necessary condition for a strong minimum.
 The simple example below demonstrates the phenomenon.

 %

\begin{example}\label{ex1}
{\rm Consider the  problem:
$$
\begin{array}{ll}
{\rm minimize}& x_2(1)\\
{\rm s.t.}& \dot x_1= u,\quad \dot x_2 = x_1\sin 2\pi u;\\
          & u\in U(t)\equiv [0,1];\\
          & x_1(0)=x_2(0)= 0.
\end{array}
$$
Then $\xb(t)\equiv (0,0)$, $\ub(t)\equiv 0$ is a weak minimum in the problem. Indeed,
$\xb_2(1)= 0$ and if $u(t)\le 1/2$ for all $t$, then both $\dot x_1(t)$ and $\dot x_2(t)$
are nonnegative for all $t$. On the other hand, it is clear that $(\xb(\cdot),\ub(\cdot))$
is not a strong minimum. Indeed,take an $\ep>0$ and $u(t)=3/4$ if $t\le\ep$ and $u(t)=0$
for $t>\ep$. Then $0>x_2(t)>-\ep$ for all $t\in (0,1])$.

Nonetheless the maximum principle and the second order condition (22)  are satisfied
but not the full secon order condition (21). Indeed, 
$H(t,x,p,u)= p_1u + p_2x_1\sin 2\pi u$, so that the adjoint system is
$$
\dot p_1= -p_2\sin 2\pi u,\qquad \dot p_2 = 0
$$
and the transversality conditions are $p_1(1)=0,\ p_2(1)= -1$,
 so that the solution
of the system corresponding to $\ub(\cdot)$ is $p_1(t)\equiv 0,\ p_2(t)\equiv -1$
and $H(t,\xb(t),p(t),u)\equiv 0$, whence the maximum principle.

We have furthermore $H_{x_1}(t,\xb(t),p(t),u)= -\sin 2\pi u$,\  $H_{x_2}(t,x,p,u)\equiv0$, \\ 
$H_u'(t,\xb(t),p(t),u)\equiv 0$, \ $H_{xx}(t,x,u)\equiv 0$, \ 
$H_{uu}(t,\xb(t),p(t),u)\equiv 0$ \\
and $H_{x_1u}(t,\xb(t),p(t),u)= -2\pi\cos 2\pi u$.

Verification of (22) is now equally simple.
The set $\Lambda$ of multipliers consists of a single element ($\la_0=1$ and $ p(t)\equiv (0,-1)$),
$T(U(t),\ub(t))\equiv \R_+$, so  the elements of the critical cone, if $\beta=0$, are
defined by the system $\dot h_1= u, \dot h_2=0$, $h_1(0)=0,\ h_2(0)=0$, $u\ge 0$.
Hence $h_1(t)\ge 0$.  $h_2(t)\equiv 0$ and $u(t)\ge 0$ for any element of the critical cone. 
The second derivative of the cost function is identical zero  and
$$
H_{(x,u)}''(t,\xb(t),p(t),\ub(t))(h(t),u(t))(h(t),u(t))= 2\pi p_2(t)\cos 2\pi\ub(t)=-2\pi h_1(t)u(t)
$$
and (22) reduces to $2\pi\int_0^1h_1(t)u(t)dt\ge 0$.

If on the other hand, elements of the critical cone associated with $w(\cdot)$ are defined by
the system $\dot h_1= u+\beta w(t) $, 
$\dot h_2=0$, $h_1(0)=0,\ h_2(0)=0$, $u\ge 0$. So if $\beta>0$, then taking $u(t)\equiv 0$ and $w(t)=3/4$ for 
$t\in [0,\ep]$ and $w(t)=0$ for $t>\ep$, we find that $\ep>h_1(t)>0$ for $t>0$ and the left-hand side of (21) reduces to
$$
\beta\int_0^TH_{x_1}(t,\xb(t),p(t),w(t))h_1(t)dt =\beta\int_0^{\ep} h_1(t)\sin(3/2)\pi dt <0.
$$
}	
\end{example}

It seems that (21) is a new type of a second order condition that has not  appeared in the literature so far. Condition (22), on the other hand,  extends a recent result
of Frankowska and Osmolovskii \cite{FO18} proved for autonomous 
problems without equality end point constraints. 
 

We furher note that, although in \cite{FO18} no analogue of (H$_9$) is explicitly stated, the
condition is automatically satisfied for the problem considered there because of the absence
of end point equality constraints. (We also mention, to avoid confusion with signs, that the $p(\cdot)$ in \
\cite{FO18} is the same as  $-p(\cdot)$ here, so  the first order condition there is the ``minimum principle''and  plus stands before the last integral in (22).)

2. Another necessary condition for a strong minimum was proved by Pales and Zeidan in \cite{PZ07}. In this result $\ub(\cdot)$ is compared with controls that may 
substantially differ from it on sets of positive measure. But these controls have very specific structure: all of them have the form $\ub(t+\theta(t))$, where $\theta(t)$
is uniformly small. Moreover, the proof of the result essentially relies on a differentiability assumption on $f$ with respect to $t$, as well as on the assumption that the control set $U(t)$ does not depend on $t$. Under these assumptions the authors,
using a modification of the method of Dubovitzkii and Milyutin, pass to another problem
in which $t$  appears as a state
variable. The first and second order necessary optimality conditions
for the strong minimum in the original problem are then obtained as the first and second order conditions necessary for a weak minimum in the new problem. 

However, as long as the first order condition is in our disposal,
a simpler construction (also involving a change of the time variable) can be used to get the second order condition of \cite{PZ07}. 
In the context of ({\bf OC2})\footnote{In \cite{PZ07} a problem with variable time interval is considered plus the inequality constraints have more general structure than here.}, with $U$ not depending of $t$,
 it is enough to consider the problem
$$
\begin{array}{rl}
{\rm minimize}& \ell_0(y(0),y(1)),\\
{\rm s.t.}&  \dfrac{dy}{d\tau}= vf(t,y,\ub(T\tau));\quad\dfrac{dt}{d\tau}=v,\quad v\ge 0;\\
&g_i(t,y)\le 0,\; i=1,\ldots,s;\\
& \ell_j(y(0),y(1))\le 0,\; j=1,\ldots,l,\\ & \ell_j(y(0),y(1))=0,\; j=l+1,\ldots,r.
\end{array}
$$
It is obvious that $\bar t(\tau)=T\tau,\; \yb(\tau)= \xb(T\tau), \vb(\tau)\equiv T$
is a local minimum in the problem. Applying (21) to this problem we get a second order condition very close to that of \cite{PZ07} with a slightly different critical cone.
We leave the details to the reader. (Just note that, although the full maximum principle for ({\bf OC2}) cannot be obtained from the first order condition for the last problem,   the adjoint equation and the transversality conditions for ({\bf OC2}) are easily recovered.)
It is also an easy matter to see that the second order condition of \cite{PZ07} is satisfied in Example \ref{ex1}.

3. There is a group of closely connected second order conditions for a strong minimum
for the case when the optimal control $\ub(\cdot)$ is piecewise continuous, in particular
for bang-bang controls (see e.g. monographs of Milyutin--Osmolovskii \cite{MO} and Osmolovskii-Maurer \cite{OM}). Our theorem that does not take specific structure of
the optimal control into account and works for arbitrary measurable controls clearly
does not cover these results. However, it is not a difficult matter to see that the
unconstrained reduction techniques developed in the proof in the next section can be
applied to optimal controls  having special structure and in particular to analyze 
conditions coming from variations of points of discontinuity of the optimal control when the latter is piecewise continuous.

\subsection{Unconstrained reduction}

It is an easy matter to see that (H$_9$)-(H$_{11}$)
together with (H$_3$) imply the hypotheses (H$_1$)-(H$_5$)
and we can use the reduction theorems for our problem.
This time we shall use Theorem \ref{redthm} (rather than Theorem \ref{red2}) applied to problem ({\bf OC2}). 
Take a small $\del>0$ as in section 3 and set $U_0(t)=B(\ub(t),\del)\cap U(t)$.
Let as before $u_1(\cdot),\ldots,u_k(\cdot)$ be a finite collection of of elements of $\cu$ (defined as in Section 3) and
$U_k(t)=U_0(t)\cup\{u_1(t),u_2(t),\ldots,u_k(t)\}$.
Then $(\xb(\cdot),\ub(\cdot))$ is a strong local minimum in the problem ({\bf OC2}$_k$)
 obtained from ({\bf OC2}) if we replace $U(t)$ by $U_k(t)$. Thanks to (H$_{12}$) we only need to consider the non-singular case.   By Theorem \ref{redthm} there is a $\la>0$ such that $(\xb(\cdot),\ub(\cdot),0,\ldots,0)$ is a local minimum of the functional
$$
\begin{array}{l}
\cj_k(x(\cdot),u(\cdot),\al_1,\ldots,\al_k)= \la\psi(x(\cdot))+\dis\sum_{i=l}^r|\ell_i(x(0),x(T))|\\
\qquad\qquad
+\dis\int_0^T\|\dot x(t)- f(t,x(t),u(t))- \sum_{i=1}^k\al_i(f(t,x(t),u_i(t))-f(t,x(t),u(t)))\|dt
\end{array}
$$
on $Z_k$ (defined in Section 3) in  the topology of  $W^{1,1}\times L^{\infty}\times \R^k$. 
It is easy to see that $\cj_k$ is a particular case of the basic model (12).

Indeed, consider the following function on $\R^{r+1}\times (C[0,T])^s\times L^1$:
$$
\begin{array}{l}
g((\xi_0,\ldots,\xi_r),(y_1(\cdot),\ldots,y_s(\cdot)),z(\cdot ))\\
\qquad\qquad\qquad=
\la \max\{\dis\max_{0\le i\le l}\xi_i,\dis\max_{1\le i\le s}\dis\max_{0\le t\le T}y_i(t)  \} + \dis\sum_{i=l+1}^r|\xi_i| + \dis\int_0^T\| w(t)\|dt.
\end{array}
$$ 
Let further $U$ be the collection of all bounded measurable selections $u(\cdot)$ of $U(\cdot)$.
Note that by (H$_{10}$) any such $u(\cdot)$ belongs to $\cu$.
Take $u_1(\cdot),\ldots,u_k(\cdot)\in U$, and let the mappings $F_i: W^{1,1}\times L^{1}\to \R^{r+1}\times (C[0,T])^s\times L^1$, $i=0,\ldots,k$ be defined as follows:
$$
\begin{array}{l}
F_0(x(\cdot),u(\cdot))=((\xi_{00},\ldots,\xi_{0r}),(y_{01}(\cdot),\ldots,y_{0s}(\cdot)), w_0(\cdot) );\\
F_i(x(\cdot),u(\cdot))=((0,\ldots,0),(0,\ldots,0),w_i(\cdot)),\; i=1,\ldots,k,
\end{array}
$$
where 
$$
\begin{array}{l}
\xi_{0j}=\ell_i(x(0),x(T));\quad y_{0i}(t)=g_i(t,x(t));\\
w_0(t)=\dot x(t)-f(t,x(t),u(t));\\
 w_i(t)= -(f(t,x(t)u_i(t))- f(t,x(t),u(t))),\ i=1,\ldots,k.
\end{array}
$$
Then $\cj_k$ is precisely $g(F_0(x(\cdot),u(\cdot))+ \sum\al_i F_i(x(\cdot),u(\cdot)) )$
and all assumptions of (H$_7$) are obviously satisfied.

A first order necessary condition for $(\xb(\cdot),\ub(\cdot),0,\ldots,0)$ to be a local minimum of $\cj_k$ can of course be obtained from Proposition \ref{pro2}. However we prefer to give an independent (and  fairly
simple) proof based on (13) with the aim to emphasize connection with the subsequent proof of the second order condition.

In terms of (12), $\yb=(\bar{\xi}_0,\ldots,\bar{\xi}_r,\yb_1(\cdot),\ldots,\yb_s(\cdot),0)$,
where $\bar{\xi}_i=\ell_i(\xb(0),\xb(T))$ and $\yb_i(t)=g_i(t,\xb(t))$, and  any $y^*\in \sd g(\yb)$ is represented by
$(\la_0,\ldots,\la_r,\mu_1,\ldots,\mu_s,p(\cdot))$, where $\la_i\ge 0$ ,
$\la_i\ell_i(\xb(0),\xb(T))=0$ for $i=0,\ldots,\ell$, \  $|\la_i|\le 1$ for $i=l+1,\ldots,r$, \
$\mu_i$  are nonnegative measures supported on $\Delta_i=\{t:\; g_i(t,\xb(t))=0  \}$, 
$\la_0+\ldots+\la_l+\mu_1([0,T])+\ldots+\mu_s([0,T])=\la$ and $\|p(t)\|\le 1$ almost everywhere.
The set $\Lambda_k$ of elements of $\sd g(\yb)$ satisfying (13), that is such that
$$
\begin{array}{l}
\dis\sum_{i=0}^r\la_i\ell_i'(\xb(0),\xb(T ))(h(0),h(T)) +\dis\sum_{i=1}^s\int_0^Tg_{ix}'(t,\xb(t))h(t)\mu_i(dt)\\
\qquad\qquad\qquad\qquad\qquad+\dis\int_0^T\lan p(t),\dot h(t)-f_x'(t,\xb(t),\ub(t)h(t)\ran dt=0
\end{array}
\eqno (23)
$$
for all $h(\cdot)\in W^{1,1}$ and
$$
\int_0^T\lan p(t),f(t,\xb(t),\ub(t))-f(t,\xb(t),u_i(t))\ran
dt\ge 0,\; i=1,\ldots,k
\eqno (24)
$$
is nonempty. Setting $h(t) = h(0)+\int_0^t\dot h(s)ds$ and changing the order of integration in the second term and the second part of the third terms of (23), and applying afterwards (23) with $h(\cdot)$ equal zero at the ends of the interval, we find that for such $h(\cdot)$
$$
\int_0^T\lan p(t)+\int_t^T (\sum_{i=1}^sg_{ix}'(s,\xb(s))d\mu_i(s) - H_x'(s,\xb(s),p(s),\ub(s))ds),\dot h(t)\ran dt= 0. 
$$
It follows that
$$
p(t)+\int_t^T (\sum_{i=1}^sg_{ix}'(s,\xb(s))d\mu_i(s) - H_x'(s,\xb(s),p(s),\ub(s))ds)= {\rm const}\quad{\rm a.e.}, 
\eqno (25)
$$
the constant obviously equal to $\lim_{t\to T}p(t)+\sum_{i=1}^sg_{ix}'(T,\xb(T))\mu_i(\{T\})$.
In turn, this implies that
$p(\cdot)$, having been corrected on a set of measure zero, if necessary, becomes a function of bounded variation. Returning back to the original form of (23) and integrating the first term in the last integral by parts, we conclude, setting $p(T)=\lim_{t\to T}p(t)$, that for any $h(\cdot)\in W^{1,1}$

$$
\begin{array}{l}
\dis\sum_{i=0}^r\la_i\lan\ell_i'(\xb(0),\xb(T ),(h(0),h(T))\ran \\
\qquad\qquad\qquad +\dis\sum_{i=1}^s \lan g_{ix}'(T,\xb(T))\mu_i(\{T\}),h(T)\ran
+\lan p(T),h(T)\ran-\lan p(0),h(0)\ran=0.
\end{array}
$$

which means that 
$$
(p(0),-(p(T) + \sum_{i=1}^s g_{ix}'(T,\xb(T))\mu_i(\{T\}))= \sum_{i=0}^r\la_i\ell_i'(\xb(0),\xb(T )) 
\eqno (26)
$$
We can now summarize.

\begin{proposition}\label{pnes} We assume (H$_9$)-(H$_{11}$). If $(\xb(\cdot),\ub(\cdot),0,\ldots,0)$ is a local minimum of $\cj_k$, then the set $\Lambda_k$ of
tuples $(\la,\ldots,\la_r,\mu_1,\ldots,\mu_s,p(\cdot))$, where $\la_i$ are numbers, $\mu_i$ are nonnegative measures supported on $\Delta_i$ and $p(\cdot)$ are functions of bounded variations, satisfying 
$$
\la_i\ge 0,\ j=0,\ldots,l,\; \la_0+\cdots+\la_l+|\la_{l+1}|+\cdots+|\la_r|+\mu_1([0,T])+\cdots+\mu_s([0,T])>0
$$
along with (24)-(26) is nonempty.
	
\end{proposition}

\subsection{Proof of Theorem \ref{th6}}

The structure of the proof does not much differ from the structure of the proof of the maximum principle in Section 4.
It is actually simpler as we do not need to consider the singular case.
We first get a second order condition for $\cj_k$ using Theorem \ref{th4}, then reformulate this condition for the problem ({\bf OC2$_k$}) and eventually for the original problem ({\bf OC2}).

So given a finite collection $\{u_1(\cdot),\ldots, u_k(\cdot)\}$ of bounded measurable selections of $U(\cdot)$, we set
$U_k(t)=U_0(t)\cup\{u_1(t),u_2(t),\ldots,u_k(t) \}$ and consider the corresponding functional $\cj_k$.

To be able to apply Theorem \ref{th4}, we need to verify that $\cj_k$ satisfies all conditions of the theorem under a suitable choice of the Banach spaces $X, \ W,\ V$ and $Y$. Set $X=W^{1,1}$,  $W=L^{\infty}$ (so that $U\subset W$), and let $Y=\R^{r+1}\times(C[0,T])^s\times W^{1,1}$. It is an easy matter to see that
the mappings $F_i: X\times W\to Y$ defined in the previous section are continuously differentiable near $(\xb,\ub)$ and twice differentiable at the point. This means that   
(H$_7$) holds for our problem.  Let further $V$ be the space of measurable $v(\cdot)$ with the norm
$$
\|v(\cdot)\|_V=\int_0^T\|f_u'(t,\xb(t),\ub(t)) v(t)\|dt.
$$
Verification of (H$_8$) now does not present any difficulty.

Thus, all we need is to verify that for a pair $(u(\cdot),v(\cdot))$
of second order feasible variations 
and  for any $\ep>0$ there are sequences of $v_m(\cdot)\in L^{\infty}$
 converging to $v(\cdot)$ in $V$ and of positive $\la_m\to 0$
such that (18) holds (with $t_m$ replaced by $\la_m$).

So let an $\ep>0$ and a pair $(u(\cdot),v(\cdot))$ of second order feasible variations be given. Take $0<\la_m\le \ep/m$,
and for any $m=1,2,\ldots$ let
$$
\Delta_m=\{t:\; \| v(t)\|\le m,\; d(\ub(t)+\la_mu(t)+\la_m^2v(t),U(t) )\le\la_m^2\ep  \}.
$$
Clearly, the measure of $\Delta_m$ goes to $T$ as $m\to \infty$.

Now for $t\in\Delta_m$ we choose a measurable $v_m(t)$  
such that $\|v(t)-v_m(t)\|\le \ep$ and
$\ub(t)+\la_mu(t)+ \la_m^2v_m(t)\in U(t)$.
For $t\not\in\Delta_m$ we define $v_m(t)$ 
also satisfying the last inclusion and such that
$\la_m^2\| v_m(t)\|=d(\ub(t)+\la_mu(t),U(t))$.  Then  $\la_m\|v_m(t)\|\le \|u(t)\|$ almost everywhere as $\ub(t)\in U(t)$. Since $u(\cdot)$ is bounded measurable, it follows 
that every $v_m(\cdot)$ is also bounded measurable, hence
belonging to $L^{\infty}=W$. Let $\Delta_m^c$ stand for $[0,T]\backslash \Delta_m$. Then by the  the second order feasibility condition:
$$
\int_{\Delta_m^c}\|\fb_u'(t)v_m(t)\|dt\le \la_m^{-2}\int_{\Delta_m^c}\|\fb_u'(t)\|d(\ub(t)+\la_mu(t),U(t))
dt\le\int_{\Delta_m^c}\xi(t)dt\to 0
$$
as $m\to \infty$ and therefore $v_m(\cdot)\to v(\cdot)$ in $V$. 

Finally, by uniform boundedness of $u(\cdot)+\la_mv_m(\cdot)$ 
$$
\begin{array}{l}
f(t,\xb(t)+\la_mh(t),\ub(t)+\la_m(t)u(t)+\la_m^2v_m(t) )\\
\qquad\qquad=\fb(t)+\la_m(\fb_x'(t)h(t)+\fb_u'(t)(u(t)+\la_mv_m(t)))\\     \qquad\qquad\qquad\qquad + (\la_m^2/2)\fb''(t)(h(t),u(t)+\la_mv_m(t))(h(t),u(t)+\la_mv_m(t)) +q_m(t)\\   
\qquad\qquad =\fb(t)+\la_m(\fb_x'(t)h(t)+\fb_u'(t)u(t))\\
\qquad\qquad\qquad\qquad + (\la_m^2/2)(\fb''(t)(h(t),u(t))(h(t),u(t))+ 2\la_m\fb_u'(t)v_m(t)) +r_m(t).
\end{array}
$$
This is exactly what we need. Indeed, $\int_0^T\| q_m(t)\|dt=o(\la^2)$ as, follows from the uniform twice differentiability assumption of (H$_{10}$). On the other hand
$$
\begin{array}{r}
\big|\fb''(t)(h(t),u(t)+\la_mv_m(t))(h(t),u(t)+\la_mv_m(t))-\fb''(t)(h(t),u(t))(h(t),u(t))     \Big|\\
\le \|\fb''(t)\|(\|h(t)+2\|u(t)\|+\la_m\|v_m(t)\|)\|\la_mv_m(t)\|,

\end{array}
$$
so that taking a  $K\ge \|\fb''(t)\|(\|h(t)+3\|u(t)\|+\ep\|)$ a.e., we see that
$\|r_m(t)-q_m(t)\|\le K\ep$ for $t\in\Delta_m$ and $\|r_m(t)-q_m(t)\|\le K$ if $t\not\in\Delta_m$.
Therefore $\int_0^T\| r_m(t)\|dt\le (K+1)T\ep$ for sufficiently large $m$. Thus $\cj_k$
does satisfy all conditions of Theorem \ref{th4} and we can apply the theorem.

Note further that
a tuple $(h(\cdot),u(\cdot),\beta_1,\ldots,\beta_k)$ with $\beta_i\ge 0$ 
belongs to ${\rm Crit} \cj_k $, if (20) holds with the last relation replaced by
$$
\begin{array}{l}
\dot h(t) = \fb_{x}'(t)h(t) + \fb_u'(t)u(t)\\
\qquad\qquad\qquad\quad\ 
+\dis\sum_{i=1}^k\beta(f(t,\xb(t),u_i(t))-f(t,\xb(t),\ub(t))).
\end{array}
$$
(cf. Remark   5.3). 

Let further a bounded measurable selection
$w(\cdot)$ of $U(\cdot)$ be given.  Set $u_1(\cdot)=w(\cdot)$, and let $\{u_2(\cdot),\ldots,u_k(\cdot)\}$ be a finite collection of bounded measurable selections of $U(\cdot)$.
Set as before 
$U_k(t)=U_0(t)\cup\{u_1(t),\ldots,u_k(t)\}$ and consider the corresponding functional $\cj_k$.
By Theorem \ref{th4} for any $(h(\cdot),u(\cdot),\beta_1,\ldots,\beta_k)$ with
$\beta_1=\beta\ge 0$ and $\beta_i=0,\ i=2,\ldots,k$ such that (20) holds, that is such that 
$(h(\cdot),u(\cdot),\beta)\in C(w(\cdot))$, we can find a tuple of multiplies $(\la_0,\ldots,\la_r,\mu_1,\ldots,\mu_s,p(\cdot))\in\Lambda_k$ such that
(21) is satisfied. But we saw in the proof of the maximum principle that
$\|lambda_k$ are compact sets, $\Lambda_{k'}\subset \Lambda_k$ if $\Lambda_{k'}$ is defined by a bigger
set of selections of $U(\cdot)$ and the intersection $\Lambda$ of all $\Lambda_k$ is nonempty.
So we can be sure that there is a $(\la_0,\ldots,p(\cdot))\in \Lambda$ such that (21) holds.

\end{document}